\documentclass[a4paper,12pt]{article}
\usepackage[utf8]{inputenc}
\usepackage[T1]{fontenc}
\usepackage{times}
\usepackage{fourier}
\usepackage{amssymb}
\usepackage{amsfonts}
\usepackage[centertags]{amsmath}
\usepackage{amsthm}
\usepackage{color}
\usepackage{newlfont}
\usepackage{cases}
\usepackage{graphicx}
\usepackage{graphics}
\usepackage{lipsum}
\usepackage{anysize}
\usepackage{float}
\usepackage{cite}
\usepackage{hyper}
\usepackage[font={small,it}]{caption}
\usepackage{multirow,authblk,lineno,bm,makecell,multicol}
\usepackage{cuted,fancyhdr}
\usepackage[dvipsnames]{xcolor}
\def\circleit#1{{\setbox0=\hbox{$\bigcirc$}\setbox1=\hbox{#1}% 
		\dimen10=\wd0 \advance\dimen10 by \wd1\divide\dimen10 by 2      
		\dimen12=\ht0 \advance\dimen12 by \dp0 
		\advance\dimen12 by-\ht1 \advance\dimen12 by-\dp1
		\divide\dimen12 by 2 \advance\dimen12 by-\dp0 \advance\dimen12 by \dp1
		\hbox to \wd0{\lower\dimen12\copy0\kern-\dimen10\copy1\hss}}}

\theoremstyle{plain}

\theoremstyle{corollary}

\theoremstyle{lemma}

\theoremstyle{proposition}
\newtheorem{proposition }{Proposition}
\theoremstyle{definition}
\newtheorem{definition  }{Definition}[section]
\theoremstyle{remark}

\theoremstyle{example}

\numberwithin{equation}{section}

\begin{document}
\author{\small{MOHAMED A. AMER}}
\title{\small{FOUNDATIONS OF TRIGONOMETRY: CONCEPTUAL AND LOGICAL
	 \\
Being an Essay Towards a Conceptual Foundations of Mathematics\thanks{%
2020 \textit{Mathematics Subject Classification}. Primary: 03A05, 00A30,
03-03, 01A20, 01A99, 33B10, 33-03. Secondary: 97E99, 97G60, 97I99, 97G30.}%
\thanks{\textit{Key words, symbols and phrases}. Conceptual soundness,
Eudoxus, Descartes, Lebesgue,  Eudoxus' theory of ratio and proportion, numbers,
measure, direct measure, indirect measure, length, area, measure of angles,
sine, geometry, analysis, $\pi $.}}}
\date{\small{TO EUDOXUS \\
c408-c355 BC\\
The Savior of Ancient Mathematics, The Founder of Modern Mathematics \\And\\DESCARTES\\1596-1650\\
The Cofounder of the Real Number System
}}

\maketitle

%\vskip0.75cm

\begin{tabular}{lll}
\qquad \qquad \qquad \qquad\qquad & \qquad \qquad  \qquad \qquad \qquad & 
[...] .In contrast, the present exposition,\\ 
&  &  in my opinion, is marred by a gross\\ 
&  &  error, not against logic but against\\ 
&  &  common sense, which is more serious.\\ 
&  & Henri Leon Lebesgue \\ 
&  & [1875-1941]
\end{tabular}

\vskip1.0cm
%\newpage
\def\baselinestretch{1.75}

\textbf{Abstract.}
	Noticing that all of the $19^{\text{th}}$, $20^{\text{th}}$ and $21^{\text{st}}$ centuries treatments of trigonometry surveyed in this article are conceptually or logically defective, it is required to seek a conceptually sound and logically correct foundations  of the subject. To this end, several questions have to be discussed: Is mathematics arbitrary? What does it have to do with nature? Reality? Applications? What is measure? Direct  measure? Indirect measure? What are ratios? Eudoxean ratios? What is their relationship to measure? What are the real numbers? What is their relationship to ratios? To measure? What are the geometric trigonometric functions? The analytic trigonometric functions? What is the relationship between them? What is the measure of an angle? what is its relationship to trigonometry?......? 
	
	After dealing with these philosophical and technical questions, a treatment of both geometric and analytic  trigonometry which would be conceptually sound and logically correct foundations of the subject is proposed.
	
	This treatment distinguishes between geometric trigonometric functions (the elements of whose domain are angles) and analytic trigonometric functions  (the elements of whose domain are real numbers), the bridge from the former to the latter is the measure of angles.
	
	The geometric  trigonometric functions are defined to be the Eudoxean ratios between appropriate straight line-segments, and the measure of an angle is defined to be its Eudoxean ratio to the radian, which is to be defined beforehand .
	
	Also, the  Eudoxean ratios were defined beforehand, and it was advocated that it is conceptually sound to consider them as the positive elements of the field of the real numbers.

\textbf{Remark.}
	Parts of this article are taken (without further notice) from the preprint [Amer [2017]] almost verbatim. Meanwhile, there are major differences between the two articles.

\vskip0.5cm

\textbf{TABLE OF CONTENTS}  
\vskip0.25cm

\textbf{I. BACKGROUND} 
\vskip0.250cm
\hskip0.75cm I.1. Question and Answers 

\hskip0.75cm I.2. Towards a Careful Analysis 

\hskip0.75cm I.3. Logical Integrity 

\hskip0.75cm I.4. Arbitrary and Accidental? 

\hskip0.75cm I.5. Direct Measure 

\hskip0.75cm I.6. Indirect Measure 

\hskip0.75cm I.7. Erroneous but Fecund Mathematics  
\vskip0.250cm
\textbf{II. FROM ANTIQUITY TO MODERNITY} 
\vskip0.250cm
\hskip0.750cm II.1. Eudoxus' Theory of Ratio and Proportion 

\hskip0.750cm II.2. Ratios and Real Numbers 

\hskip1.0cm II.2.1. The Eudoxus-Descartes Field of Real Numbers (The E-D Field)

\hskip1.5cm II.2.1.1. Ratios 

\hskip1.5cm II.2.1.2. Operations

\hskip1.5cm II.2.1.3. Descartes, Hilbert and Tarski-Zero and the Negatives

\hskip1.5cm II.2.1.4. Continuity

\hskip0.750cm II.3. Angles 

\hskip0.750cm II.4. Measure of Angles 

\hskip0.750cm II.5. From Geometry to Analysis  

\hskip0.750cm II.6. Alternative Approaches 
\vskip0.250cm
\textbf{III. CONCLUDING REMARKS } 
\vskip0.250cm
\textbf{ ACKNOWLEDGMENTS}
\vskip0.25cm
\textbf{APPENDIX}
\vskip0.25cm
\textbf{BIBLIOGRAPHY}

\newpage

\textbf{I. BACKGROUND}
\vskip0.25cm
\textbf{I.1 Question and Answers.} 
That the foundations of trigonometry is problematic may be illustrated by the observation  raised by Hardy
[1967, p.    316] over a century ago: ``The whole difficulty lies in the
question, \textit{what is the } $x$\textit{\ which occurs in }$\cos x$ and $\sin x$?'' 
"To answer this question," continues Hardy [1967, pp.     316-7]:  
\begin{center}
\begin{tabular}{p{4.5in}}   
	``we must define the measure of an angle, and we are now in a
	position to do so. The most natural definition would be this: suppose that $%
	AP$ is an arc of a circle whose centre is $O$ and whose radius is unity, so
	that $OA=OP=1$. Then $x$, the measure of the angle [$AOP$], is \textit{the
		length of the arc} $AP$. This is, in substance, the definition adopted in
	the text-books, in the accounts which they give of the theory of `circular
	measure'. It has however, for our present purpose, a fatal defect; for we
	have not proved that the arc of a curve, even of a circle, possesses a
	length. The notion of the length of a curve is capable of precise
	mathematical analysis just as much as that of an area; but the analysis,
	although of the same general character as that of the preceding sections, is
	decidedly more difficult, and it is impossible that we should give any
	general treatment of the subject here.
	
	We must therefore found our definition on the notion not of length but of 
	\textit{area}. We define the measure of the angle $AOP$ as \textit{twice the
		area of the sector }$AOP$\textit{\ of the unit circle.}''\\
	\end{tabular}
\end{center}

Applying the area approach to measuring angles, hardy [1967, p.    317] practically shows that for an (acute) angle $\theta$, the measure of $\theta$ is equal to $\int_{0}^{m}\frac{dt}{1+t^2}$, where $m$ is the (geometric) tangent of  $\theta$; while Morrey [1962, pp.     214-8] shows that the measure of the same angle is $x\sqrt{1-x^2}+2\int_{x}^{1}\sqrt{1-t^2}\, dt$, where $x$ is the (geometric) cosine of $\theta$.

The arc-length approach to measuring angles was pursued from
different points of view, in Anton et al. [2009], Gearhart and Shultz
[1990], Gillman [1991], Hass et al. [2020], Hughes-Hallett et al. [2005], Protter and Morrey
[1963], Richman [1993], Rose [1991], and Ungar [1986].

Inspired by the arc-length approach, the treatment in Eberlein [1966] is
based on the theory of complex valued functions on real numbers. 

According to Hardy [1967, p.    433], by defining the measure of an angle (as the length of the corresponding are or twice the area of the corresponding sector [Hardy [1967]  pp.  316-7]),  trigonometry will rest on a secure foundation via translating the geometrical language employed by the ordinary textbooks into the language of analysis. He calls this approach geometrical and considers it be the most natural.

Besides, Hardy [1967, pp.     433-4] mentions four more approaches to answer his
question:
\begin{itemize}
	\item[(i)] Power series.
	\item[(ii)]  Infinite products.
	\item[(iii)]  Complex analysis. This approach was pursued, from two different points
	of view, in Davis [2003] where arc length was made use of (but in an
	avoidable way) and in Bartle and Tulcea [1968] where an ad hoc definition of
	measuring angles is introduced and made use of.
\item[(iv)]   Definition of the inverse functions by integrals. Hardy [1967, p.    434] made use of this approach to define $\tan^{-1}$ and Morrey [1962, pp.     214-8] made use of it to define $\cos^{-1}$.
\end{itemize}

  The matter in the nineteenth century was different. In order to define the measure of angles, Hobson [1918, p.   3; $1^{st}\ ed$. 1891] says ``we must decide upon a unit angle, which may be any arbitrarily chosen angle of fixed magnitude;  then all other angles will be measured numerically by the rations they bear to this unit angle.''

  It is to be noticed that the ratios here are between angles, which are geometric object, not numbers; but (according to Hobson [1918, p.   3]) the ratios are numbers. De Morgan [1837, pp.    12-14] makes use of essentiality the same definition, without explicitly stating it.

  The unit angle chosen for theoretical purposes by each of Hobson [1918, p.   4] and De Morgan [1837, p.   13] is the radian. While for practical purposes (think, e.g. , of the latitudes and longitudes; this is plane, not spherical, trigonometry) the sexagesimal system based on the degree as the unit is adopted [De Morgen [1837]p.   13; and Hobson [1918] pp.    3,4].

Besides, some books on the foundations of geometry, presupposing real
numbers, take an axiomatic approach to define the measure of angles [Moise
1964] p.    74, Borsuk and Szmielew 1960] p.    172, and Greenberg 1980] p.    98].

We shall come back to this in I.3, I.6 and II.1, II.4 below.

\vskip0.5cm
\textbf{I.2 Towards a Careful Analysis.} 
The
 foundations of the theory of trigonometric functions is not as simple as it
 is generally supposed. It still needs to be properly and carefully analyzed
 [cf. Hardy [1967] p.    432].
 
 Proper analysis must deal not only with the logical aspects, but also with
 the \textit{conceptual} aspects (which include Lebesgue's ``common sense'').
 We shall come back to this in I.4 below.
 
\vskip0.5cm
\textbf{I.3 Logical Integrity.} 
Most common approaches to the
theory of trigonometric functions are based on the length of a circular arc
or the area of a circular sector.

From the logical point of view, these approaches involve two problems, each
potentially entails circularity:
\begin{enumerate}
	\item[(i)] The trigonometric functions are defined in terms of arc length or sector
	area, while these entities are, possibly, defined or evaluated in terms of
	the trigonometric functions [see I.6 below].
	
	\item[(ii)] The proof of the celebrated limit, $\displaystyle \lim_{x\rightarrow 0}
	(\sin x)/x=1$, possibly, depends on properties of arc length or sector area
	which, in their turn, depend on properties of the trigonometric functions
	which are based on this limit [cf. Richman [1993], Rose [1991], and Ungar
	[1986]].
\end{enumerate}

Whether there is actual circularity, no circularity or the situation is
obscure, can be seen only through the detailed analysis of the treatment.

Though the two aforementioned problems are  avoided via the approaches of Davis [2003], Eberlein [1966], Hardy [1967] and Morrey [1962],
they still persist, without even being pointed out to, in modern calculus books
[cf. Anton et al. [2009], Hass et al. [2020], Hughes-Hallett et al. [2005], Protter and Morrey
[1963]].

The approaches of De Morgen [1837] and Hobson [1918] suffer from two serious defects. The \textit{radian} is ill-defined [cf. De Morgan [1837, p.   13, Hobson [1918, p.   4] and the \textit{ratio} is not defined; it is made use of along the lines of Euclid [1956] which is problematic as was noticed by De Morgan [1836, P62].

\vskip0.5cm 
\textbf{I.4 Arbitrary and Accidental?} 
Is mathematics arbitrary
and its effectiveness in natural sciences accidental?

That mathematics is arbitrary is implied by Hardy, not only when he
arbitrarily defined the measure of an angle [see I.1 above] but also when he
said that real mathematics ``must be justified as art if it can be justified
at all'' [Hardy [1969] p.   139, also cf. Griffiths [2000] p.    4].

That the effectiveness of mathematics is accidental is implied by E. p.   
Wigner when, under the title ``The Unreasonable Effectiveness of Mathematics
in the Natural Sciences'' he wrote ``The miracle of the appropriateness of
the language of mathematics for the formulation of the laws of physics is a
wonderful gift which we neither understand nor deserve.'', as was noticed by
Kac and Ulam [1971, p.    166]. Kac and Ulam agree, saying that this
effectiveness ``remains perhaps a philosophical mystery'' [Kac and Ulam
[1971] p.    12] and that it ``may be philosophically puzzling'' [Kac and Ulam
[1971] p.    191]. Yet, they give [Kac and Ulam [1971] p.    191] the following
clue for solving this puzzle ``One reason for it is the necessary condition
that measurements, and thus much of the discussion in physics and astronomy,
can be reduced to operations with numbers.''

Kac and Ulam are, in a sense, following Lebesgue in emphasizing the
importance of measure of quantities, of which he says [Lebesgue [1966] pp.    
10-11] ``There is no more fundamental subject than this. Measure is the
starting point of all mathematical applications [...] it is usually supposed
that geometry originated in the measure of [angles, lengths,] areas and
volumes. Furthermore, measure \textit{provides} [emphasis added] us with
numbers, the very subject of analysis.''

Unlike Wigner and Kac and Ulam, Lebesgue is \textit{not} puzzled. For him
``[...] those whom we have to thank for such [practically effective]
abstract considerations have been able to think in abstractions and at the
same time to perform useful work precisely because they had a particularly
acute sense of reality. It is this sense of reality that we must strive to
waken in the young [and, a fortiori, the grownups].'' [Lebesgue [1966] p.    11].

Concerning the role of conceptual considerations in developing abstract
mathematics and its foundations, Marquis [2013] may be consulted.

As a mathematician, Feferman distinguishes between \textit{structural}
and \textit{foundational} axioms of mathematics (or parts thereof) holding
that both sorts of axioms are \textit{not arbitrary} [cf. Feferman [1999] p.   
100 and Feferman et al. [2000] pp.     403, 417]. Besides, he calls for a
``philosophy grounded in inter-subjective \textit{human} conceptions [...]
to explain the apparent \textit{objectivity} [emphasis added] of
mathematics.'' [Feferman [1999] p.    110].

Though on different philosophical grounds, Feferman was anticipated by
G\"{o}del [1964b, p.    264] who asserts that Mahlo ``axioms show clearly, not
only that the axiomatic system of set theory as used today is incomplete,
but also that it can be supplemented \textit{without arbitrariness}
[emphasis added] by new axioms which only unfold the content of the concept
of set explained above.''

Concerning the rigorous unfolding of informal concepts Lavers [2009] may be
consulted. And concerning the relationship between intuitive mathematical
concepts and the corresponding formal systems Baldwin [2013] may be
consulted.

Bertrand Russell went even further than G\"{o}del and Feferman when he once
said ``Logic [hence, mathematics] is concerned with the real world just as
truly as zoology, though with its more abstract and general features,'' as
was brought to my attention by G\"{o}del [1964a, pp.     212-3]. Lukasiewicz
[1998, pp.     205-7] provides logico-philosophical grounds for this view.

Lebesgue [1966, p.    101] distinguishes between ``terminology'' which is
arbitrary, or free, to use his terminology, and ``definitions'' which are
not, or ``At the very least, some definitions, those that are meant to make practical concepts more precise, are not free.'' Frege was more explicit when he said 
[Frege [1968] pp.     $107^e-8^e$]  ``[E]ven the mathematician cannot create things at will, any more than the geographer can; he too can only discover what is there and give it a name.''

Moschovakis [1995, p.    753] seconds, asserting that for a positive
[integrable] function $f:[a,b]\rightarrow \Bbb{R}$ we do not mean the
equation:
\vspace{10pt}

\begin{equation}
\displaystyle \int_a^b f(x) \, dx= \text{the \, area\, under \, $f$ \, and \,
	above \, the \, $x$-axis}\tag{*}
\end{equation}

\vspace{10pt} 
\hspace{-15pt}``as a definition of area, for, if we did, we \textit{could not} [emphasis
added] use it to compute real areas; and we cannot prove it conclusively$
^{1} $ without a separate, precise definition of ``area''- which we might
give, of course, but whose connection with actual ``physical'' or
``geometrical'' area we would then need to justify.'' We shall come back to
this and to footnote ''1'' in  II. 1  below.

  Decades before Moschovakis, Frege [2016, vol. II, p.    100 ] said ``Now it is applicability alone which elevates arithmetic [a mathematical discipline] above a game to the rank of a science. Applicability thus necessarily belongs to it. Is it appropriate, then, to exclude from arithmetic [a mathematical discipline] what it needs to be a science?''

  Regarding applicability, Frege makes two observations concerning numbers, which are relevant to trigonometry too. The first [Frege [2016] vol. II, P157] is ``At the same time, however, we avoid the emerging problems of the latter [formal] approaches, that either measurement does not feature at all, or that it feature without any internal connection grounded in the nature of the number itself, but is mearlly tacked on externally, {\dots}'' The second [Frege [2016] vol. II, p.    85] is `` If we look more closely, we realise that a number-sign cannot on its own designate a length, [angle], force, and so on, but only in connection with the designation of a measure, a unit [of the same kind as the magnitude to be measured] such as meter [degree, radian], gramme, and so on.''

  For further discussion concerning the relationship of mathematics to
  reality, Amer [1980] and Maddy [2008] may be consulted. At any rate, this
  relationship should be a subject not only of philosophical reflection, but
  also of factual and historical scrutiny.
  
  For a general view of the problem of providing foundations Sher [2013] may
  be consulted.

  How would trigonometry be looked at in the light of the above discussion?

  In all of the ${20}^{th}$ and ${21}^{st}$ centuries approaches mentioned above (with exception of what Hardy [1967, p.   433] calls geometrical) measurement does not feature at all. In each of the two branches of the geometrical approach, it features without any internal connection  grounded in the nature of trigonometry itself, but is merely tacked on externally. In one branch of this approach  [Hardy [1967] pp.     316-7] the measure of an angle is defined to be twice the area of the corresponding sector of the unit circle; why multiplying by the factor 2, not by any other factor? In the other branch the measure is defined to be the length of the corresponding arc of the unit circle; why is  it not multiplied by some factor?    

  To start with, what is the relationship between the measure of sectors and arcs on the one hand, and the measure of angles on the other hand? In all of the ${20}^{th}$ and ${21}^{st}$ centuries approaches mentioned above (geometrical or not) there is no mention of a unit angle. As was mentioned before, the matter in the ${19}^{th}$ century was different. As a matter of fact, it is possible -under certain conditions- to measure a magnitude by measuring a corresponding magnitude of a different kind [De Morgan [1837] pp.     14-5]. We shall come back to this in II below.

  The term ``measure'' or one of its derivatives was made use of above several times; in what sense? In fact there are in the literature three sorts of measure: direct (e.g. to measure a straight line-segment by a straight line-segment), indirect (e.g. to measure a curved arc by a straight line-segment) and abstract (e.g. measure on a Boolean algebra). Direct and indirect measure are-respectively- the topics of I.5 and I.6 below.

  To conclude this subpart it is worth noting that -concerning trigonometry- the matter is not as arbitrary as it may seem. The different ${20}^{th}$ and ${21}^{st}$ options mentioned above are just different approaches to define and deal with the trigonometric functions inherited from the ${19}^{th}$ century in a way that would-be logically acceptable and, at the same time, avoids the Greek ratios which were involved in the ${19}^{th}$ century treatment. 

  \vskip0.5cm
  
   \textbf{I.5 Direct Measure.}  
The direct sense of the measure of a
   magnitude is the ratio of this magnitude to a standard (measuring) unit of
   the same kind. So, in this sense, measure is a special ratio.
   
   The founder of the theory of ratio and proportion is Eudoxus. He developed
   this theory to deal with all kinds of magnitudes, not only those of (Euclidean)
   geometry, and with the incommensurables as well as the commensurables [cf.
   Euclid [1956] vol. 2, p.    112], making use only of whole numbers. Dealing
   with the incommensurables practically produces irrational numbers, e.g. $\pi 
   $. We shall come back to this in  II. 1   and  II. 2  below.
   
   Unlike Eudoxus, Lebesgue starts with measure, then generalizes to ratios.
   Having at his disposal the decimal system of numeration (which he highly
   esteems, Lebesgue [1966] p.    18]) Lebesgue [1966, pp.     19-20] defines the
   measure of straight line-segments as follows.
   
   Let a straight line-segment $u$ be the standard measuring unit, and let $b$
   be a generic straight line-segment. Lay off $u$ on $b$ several times in the
   obvious way to determine the unique natural number $n_{o}$ such that $
   n_{o}u\leq b<(n_{o}+1)u$. This process may be performed, and the uniqueness
   of $n_{o}$ proved, in Euclidean geometry. Its termination is guaranteed by
   the Archimedean principle (which Archimedes himself attributes to Eudoxus
   [Heath [1963] p.    193 and Euclid [1956] vol. 3, pp.     15-6]).
   
   Assume $n_{o}u<b$. By well known Euclidean methods divide $u$ into $10$
   equal (i.e. congruent) parts and determine the unique natural number $n_{1}$
   ($\leq 9$) such that $(n_{o}+n_{1}10^{-1})u\leq b<(n_{o}+(n_{1}+1)10^{-1})u$%
   . Iterate this procedure whenever possible. If the iteration terminates at
   stage $m$, the length of $b$ in terms of $u$ is, by definition, 
   $n_0.n_1\cdots n_{m}$. 
   Otherwise the perpetual procedure
   will give rise to an infinite decimal, which is \textit{itself} (by \textit{%
   	definition}), at the same time, both a \textit{real number} and the \textit{%
   	length} of $b$ in terms of $u$. In this sense direct measures, here lengths
   of straight line-segments, \textit{produce} real numbers, while making use
   only of natural numbers.

  Lebesgue [1966, p.   20] asks, ``\textit{whether every unlimited sequence of numbers [digits] extending to the right and containing a decimal point}[and not eventually 9's] \textit{is a number,}'' that is, whether such a sequence necessarily originates from comparison of a segment b with the unit length $u$. He [1966, pp.     20-1] answers his question positively making use of ``the axioms of geometry (either [explicitly] stated as [implicitly] understood).'' Specifically, the existence of b is proved making use of the axiom of continuity, and the uniqueness is proved making use of the Archimedean axiom. The axiom of continuity adopted by Lebesgue [1966, pp.     20-1] states that the intersection of a sequence of nested non-empty closed intervals is non-empty. 

  This together with the Archimedean axiom, which is also adopted by Lebesgue [1966, p.    21], are equivalent to Dedekind's axiom of continuity.

  Moreover, Lebesgue  [1966, p.    22-4] defines operations (addition, multiplication, {\dots}.) on his numbers. His definitions are essentially geometric; they parallel the corresponding definitions gives by Descartes [1954, pp.     2,5].

  Furthermore, Lebesgue [1966, pp.     35-8] generalizes lengths to ratios via replacing the
  unit $u$ by a generic segment $c$, as the consequent of the ratio.
  
  The decimal system may be generalized to what may be called ``positional
  system'' by replacing $10$ by any natural number $k\geq 2$.
  
  For angles, the value of $k$ which facilitates geometric constructions most
  is $2$. Replacing $10$ by $2$ and repeating everything else almost verbatim,
  measures of angles (in terms of a standard unit angle) and ratios among them
  may be discussed, giving rise to similar results as above [cf. Lebesgue
  [1966] p.   38 and  II. 1  below].
  
  Not only angles and straight line-segments are measured in this way, but
  also time intervals, masses, ... [cf. Lebesgue [1966] p.    39].
  
  Theoretical measure is not quite the same as measure in everyday life, or
  even in scientific laboratories. As the former is meant to make the latter
  precise, the former is an abstraction and idealization of the latter, and in
  this sense it is not arbitrary, but it is essentially dictated by our
  understanding of the actual facts of the real world [cf. I.4 above].
  
  For comparison of the above treatment with that of Eudoxus see  II. 1, II. 2  below.

\vskip0.5cm
\textbf{I.6 Indirect Measure.} 
 What makes direct measure, and
 ratios in general, possible is that the two terms of the ratio are capable,
 when multiplied, of exceeding one another [see II. 1 below].
 
 Let $c$ be a circular arc, and $u$ be a straight line-segment. Within
 Euclidean geometry, none of $c$ and $u$ is capable of exceeding the other,
 no matter how (finitely) many times it is multiplied. So, accepting the
 Euclidean notion of multiplicity [see  II. 1  below] the measure (here, the
 ``length'') of $c$ in terms of $u$ can but be indirect. This concept should
 be carefully analyzed if the arc-length approach to measuring angles is to
 be made precise.
 
 Indirect measure is probably at least as ancient as the XII dynasty of the
 ancient Egyptian middle kingdom (c1900 B.C.). It is believed that it was
 shown, via an example (which may be generalized) and approximate
 calculations, that the area of the surface of a hemisphere is (exactly)
 twice as much as that of the corresponding greater circle [Vafea [1998] pp.    
 31-3].
 
 The first known treatment of indirect measure after the incommensurability
 crisis is that of Archimedes. In his \textit{Measurement of a circle},
 Archimedes shows that (i) ``the area of a circle is equal to that of a right
 angled triangle having for perpendicular the radius of the circle and for
 base its circumference.'' And (ii) ``the ratio of the circumference of any
 circle to its diameter is $<3\frac{1}{7}$ but $>3\frac{10}{71}$.'' [Heath
 [1963] p.    305].
 
 Within the framework of Euclidean geometry [Euclid 1956], each of the
 statements (i) and (ii) does not make sense, simply because the
 circumference of a circle is not straight. It is surprising that Heath
 [1963] did not raise this point, while Al-Khwarizmi [1939, pp.     55-6] and
 Descartes [1954, p.    91] practically raised it.
 
 Stein [1990, p.    180] sheds more light on this issue: ``In phys: VII iv 248a$%
 _{18}$-b$_{7}$ he [Aristotle] contends that a circular arc cannot be greater
 or smaller than a [straight] line-segment, on the grounds that if it could
 be greater or smaller, it could also be equal.''  For other views see [Stein
 [1990] footnote 48, pp.     208-9].
 
 Under the urge of practical applications, and guided by physical insight and
 geometric intuition, the work of Archimedes was incorporated into mainstream
 mathematics some way or the other.
 
 In modern mathematics this is done through developing classical analysis and
 Euclidean geometry (with continuity). It is well known that these theories
 can be formalized in Zermelo-Fraenkel set theory with choice. However,
 almost all of the scientifically applicable portions of them may be
 formalized in systems which are as weak as Peano Arithmetic [cf. Feferman
 [1999] p.    109].
 
 Concerning rectifiable curves, the common practice, partially following
 Archimedes, is to adopt a treatment which concentrates only on inscribed
 polygons. Even an advanced book like [Buck [1965] p.    321] does this.
 Lebesgue [1966, p.    97] describes this as being ``hypnotized by the word
 'inscribed'.''
 
 The theory of rectifiable curves thus treated, is susceptible to a paradox
 (similar to Schwarz paradox concerning the area of curved surfaces) which
 entails that $\pi =2$ [Lebesgue [1966] p.    98].
 
 Inspired by his physical insight, Lebesgue [1966, pp.     98-106, 117-9]
 provides for plane curves two alternatives which are free from this paradox.
 The first [Lebesgue [1966] pp.     103-5] depends on measuring angles, which
 makes the arc-length definition of the measure of an angle logically
 circular. The second [Lebesgue [1966] pp.     117-9] depends on measuring areas,
 which makes the arc-length definition of measuring angles less direct, hence
 conceptually more questionable, than the area definition. We shall come back
 to this in  II. 1, II. 4  below.
 
 Accepting the traditional arc-length definition of measuring angles, and the
 traditional definition of arc-length, it is worthwhile to note that the
 celebrated limit,
 $\displaystyle \lim_{\theta \rightarrow 0} (\sin \theta
 )/\theta =1$, or rather $\displaystyle \lim_{\theta \rightarrow 0} 2\theta
 /(2\sin \theta )=1$, is a special case of the evident: 
 $$ \displaystyle \lim_{h \rightarrow 0} \int_{x}^{x+h} \sqrt{
 	1+f^{^{\prime}2}(u)}du/\sqrt{h^{2}+(f(x+h)-f(x))^{2}}=1,$$ 
 where $f$ is continuously differentiable [cf. Hardy [1967] p.    270].
  
\vskip0.5cm
\textbf{I.7 Erroneous but Fecund Mathematics.} 
   The aforementioned
  work of Archimedes and the current treatment of the trigonometric functions
  are examples of (logically or conceptually) erroneous but fecund mathematics.

  The most widely spread example of erroneous but fecund mathematics is the pre-Cauchy calculus. The following passage shows how calculus was founded by one of the most prominent mathematicians of that period [Euier [2000, pp.     51-2],  ``If we accept the notation used in the analysis of the infinite,
  then $dx$ indicates a quantity that is infinitely small, so that both $dx=0$
  and $adx=0$, where $a$ is any finite quantity. Despite this, the geometric
  ratio $adx:dx$ is finite, namely $a:1$. For this reason the two infinitely
  small quantities $dx$ and $adx$, both being \textit{equal} [emphasis added]
  to $0$, \textit{cannot} [emphasis added] be confused [\textit{considered
  	equal}] when we consider their ratio. In a similar way, we will deal with
  infinitely small quantities $dx$ and $dy$. \textit{Although these are both
  	equal to }$0$, \textit{still their ratio is not that of equals} [emphasis
  added]. Indeed the whole force of differential calculus is concerned with
  the investigation of the ratios of any two infinitely small quantities of
  this kind. The application of these ratios at first sight might seem to be
  minimal. Nevertheless, it turns out to be very great, which becomes clearer
  with each passing day.''

  Although multiplying surd roots of natural numbers has been known since long ago [cf. Al-Khwarizmi [1939, p.   32], Dedekind  
-after defining the usual operations on his real numbers- says [Dedekind [1901, p.   11],
  ``[A]nd in this way we arrive at real proofs of theorems (as, e.g., $\sqrt{2}. \sqrt{3}=\sqrt{6}$), which to the best of my knowledge have never been established before.'' On what basis was Al-khwarizmi making use of these theorems?
 Are there real proofs and unreal proofs?...?

  Correct mathematics may be sterile. The first attempt to algebraize logic,
  due to Leibniz [cf. Kneale and Kneale [1966] pp.     320-45 and Amer [2022] pp.     83-4] is a good example.
  The second attempt, due to Boole [1948], is well known to be erroneous but
  fecund.

  A more important example is Frege's system of logic, which was proved by Russell to be self-contradictory, while it proved to be extremely fecund. 

  Erroneous (or defective in some way or the other) but fecund mathematics has played an essential role in the historical development of mathematics. I do not think that this is satisfactorily taken into account by the current philosophy of mathematics.  
  
  The subject needs to be further investigated through concerted endeavors of
  mathematicians, philosophers of mathematics and historians of mathematics.
  
  \vskip1.0cm

\textbf{II. FROM ANTIQUITY TO MODERNITY}
\vskip0.50cm
\textbf{II.1 Eudoxus Theory of Ratio and Proportion.} 
That the
theory of ratio and proportion expounded in Book V of Euclid's Elements
[Euclid [1956] vol. 2, pp.     112-86] is due to Eudoxus, is beyond reasonable
doubt [Heath [1963] p.    190]. Nevertheless, the actual arrangement and
sequence of book V is not attributed to Eudoxus, but to Euclid [Heath [1963]
p.    224]. 

It should be emphasized that it was understood from the very beginning that
Eudoxus' theory of ratio and proportion (henceforth ETRP) is the foundation
not only of geometry, arithmetic (see below) and music (harmonics) but also of \textit{%
	all} mathematical sciences [Euclid [1956] vol. 2, pp.     112-3]. At Eudoxus' time
these sciences used to include, in addition, mechanics, astronomy and optics
[cf. Amer [1980] p.    569].

Moreover, it is to be noted that Euclid distinguishes between two sets of
indemonstrable principles: postulates [Euclid [1956] vol. 1, pp.     154-5] and
common notions [Euclid [1956] vol. 1, p.    155]. Postulates are peculiar to
geometry, in contrast, common notions, which are instrumental in ETRP, are
common to all demonstrative sciences [cf. Euclid [1956] vol.1, p.    221].

The key definitions of ETRP are definitions V.3 of Euclid - V.7 of Euclid
(this is an abbreviation of ``definitions 3-7 of book V of Euclid [1956]'';
in the sequel similar abbreviations will be made use of in similar
situations, without further notice).
 
{DEFINITION V.3 of Euclid}
	[1956, vol. 2, p.    114]. A \textbf{ratio}  [ratioal relation] is a sort
	of relation in respect to size between two magnitudes of the same kind.

{DEFINITION V.4 of Euclid}
 [1956, vol. 2, p.    114]. Magnitudes are said to 
	\textbf{have a ratio} [ratioal relation] to one another which are capable, when multiplied, of
	exceeding one another.

Magnitudes do not have to be geometric. They may be parameters of physical
objects: masses, electric charges, time intervals, ... .

If definition V.3 of Euclid is understood to the effect that being of the
same kind is a necessary and sufficient condition for two magnitudes to have
a ratioal relation; definition V.4 of Euclid may be considered as a definition of
``being of the same kind'', or simply, as a definition of ``Kind''. Though
Heath says that this view (which is adopted here) is accepted by De Morgan,
he sees the matter differently [Euclid [1956] vol. 2, p.    120].

In definition V.4 of Euclid ``multiplied'' means multiplied by a positive
integer (or a non-zero finite cardinal): to multiply a magnitude by 1, is to leave it as it is, to multiply
it by $n+1$ is to add it to itself $n$ times [cf. Euclid [1956] vol. 2, pp.    
138-9, also see below]. This may be made precise via recursion.

The term ``added'' occurs in the formulation of common notion 2 [Euclid
[1956] vol. 1, p.    155] which is made use of from the very beginning, e.g.,
in Pythagoras theorem [Euclid [1956] vol. 1, pp.     349-50]. In that usage ``addition'' means essentially-disjoint (i.e., only boundary points may
be in common) union [cf. Euclid [1956] vol. 1, pp.     349-50 and  I.5  above,
also see below].

To appreciate the subtleness of definition V.4 of Euclid, the following may
be easily seen in the light of the above discussion. There is a ratioal relation 
between a side of a square and its diameter and there is a ratioal  relation between a
circle and the square on its radius. Also there is a ratioal relation between any two
circular arcs of any two congruent circles. In contrast, at least within
Euclidean geometry, there is no ratioal relation between any circular arc and any
straight line-segment [cf. I.6 above]. Also arcs of non-congruent circles do
not have ratioal relations to one another, nor do elliptic arcs have a ratioal relation to one
another, even if they belong to congruent ellipses (except in some obvious
special cases). Whenever there is a ratioal relation, it is direct, no rectification or
integration is involved.

Calling it a ``species of quantity'' Stein [1990, p.    167] formalized
Eudoxus' notion of a kind of magnitude (henceforth a kind) [cf. De Morgan [1836] pp.     2-6 and Avigad et
al. [2009] p.    722]. Slightly modifying Stein's formalization and requiring
the Archimedean principle to hold, an ordered quadruple $\mathbf{K}=<K,+,<,%
\stackrel{.}{-}>$, where $+$, $<$, $\stackrel{.}{-}$ are, respectively, a binary operation, a binary relation, and a ternary relation on $K$,   is said to be a kind of magnitude or simply, a kind if it satisfies:
\begin{enumerate}
	\item[(i)] $<K,+>$ is a commutative semigroup.   
	\item[(ii)] $x<y\leftrightarrow \exists z(y=x+z).$ Hence $<$ is transitive.
	\item[(iii)] $<$ is irreflexive and trichotomous.
	\item[(iv)] $\stackrel{.}{-}(x,y,z)\leftrightarrow x=y+z$. In this case it may be
	written that $x\stackrel{.}{-}y=z$,
	 for uniqueness see below.
	\item[(v)] $<$ is Archimedean.
\end{enumerate}

This entails that $\mathbf{K}$ is an Archimedean strictly linearly ordered
cancellative commutative semigroup with partial subtraction. $K$ is said to
be the universe of $\mathbf{K}$, and is denoted also by ``$|\mathbf{K}|$''. As usual, the converse of $<$ is denoted by $">"$. Similar terminology will henceforth be made use of without further notice.

A simple example of a kind is $\Bbb{N}^{*}$ (the non-zero natural numbers)
with the usual addition ($<$ , $\stackrel{.}{-}$ are definable).
Nevertheless Euclid [ 1956] only implicitly considered it to be so [cf. Prop.    X.5 of Euclid [1956, vol. 3, p.    24]]. Instead of the general
definition of proportion (definition V.6 of Euclid [1956, vol. 2, p.    114]), $%
\Bbb{N}^{*}$ was given a special one (definition VII.20 of Euclid [1956,
vol. 2, p.    278]) which is equivalent, in this case, to the general
definition but easier to apply.

Neither the cardinals with cardinal addition, nor the ordinals with ordinal
addition, nor the positive reals with the usual multiplication can be
expanded to become a kind.

However, keeping the intuitive notion of the number (of the elements) of a
set in mind and guided by Euclid's common notions (henceforth cns, and cn
for the singular) [Euclid [1956] vol. 1, p.    155], a set theoretic kind which
is isomorphic to the one given above may be obtained as follows. In
accordance with cn 4 ``Things which coincide with [congruent to, bijective
to, ...] one another are equal [equivalent] to one another.'',  equivalence
classes of sets under bijection are considered (proper classes may be replaced by sets via the axiom of regularity). This is the only
conceptually sound choice as long as applying numbers to sets is the
objective [cf. G\"{o}del [1964b] pp.     258-9]. And in accordance with cn 5
``The whole is greater than the part.'' only sets which are not bijective to
any of its respective proper subsets are taken. Addition is readily
definable via disjoint unions. To simplify, consider the set of all finite subsets of an infinite set. Take the universe of the kind to be the quotient set under bijection, addition is readily definable via disjoint unions.

In the forthcoming applications it is needed to deal with five different kinds:
straight line-segments (henceforth slss, and sls for the singular), arcs of
congruent circles, angles, polygons and, to generalize, plane regions
bounded by slss and circular arcs.

These kinds are explications of corresponding implicit kinds in Euclid
[1956]. In all cases the explications as well as the corresponding implicit
notions are not arbitrary inasmuch as they are generalizing the
corresponding practical concepts (ratios, hence direct measures) and making
them precise [cf. Lebesgue [1966] pp.     45, 65, 101 and Borceux [2014] pp.    
vii, 45, 305]. This abstract formalization which is grounded in practical
dealing with reality [cf. Lebesgue [1966] p.    11, see I.4 above] probably
provides the justification Moschovakis [1995, p.    753, see I.4 above] is
seeking to connect mathematics with the \textit{real} and the \textit{actual}%
, which is an essential issue regarding conceptual soundness.

To see how this would go, consider footnote ``1'' of Moschovakis [1995, p.   
753, see I.4 above]. It reads, in part, as follows ``We can almost prove
(*), by making some natural assumptions about ``area'' (that we know it for
rectangles and that it is an additive set function when defined), and then
showing that there is a unique way to assign area to ``nice'' sets subject
to these assumptions.''

The problem is that the assumption that we know the area of rectangles needs
justification. And, which is more basic: why is it taken for granted that
the area function is real valued?

In his first exposition, Lebesgue [1966, pp.     42-5] deals with rectangles as
a subset of a set of regions which includes circles too. The area of
rectangles and that of the other regions are treated on an equal footing, pre-assuming the real numbers and 
making use of \textit{limits}. Besides, Lebesgue [1966, pp.     51-9] presents
other expositions dealing with the area of polygons only. In each of them
the well known formula for the area of a triangle is \textit{preassumed}.
Likewise, Moise [1964, p.    168] begins building his area function for
polygons by taking this well known formula as a \textit{definition} of the
area of a triangle\textit{. }Earlier he [Moise [1964] pp.     154-5] adopts the
well known formula for the area of a rectangle as one of the \textit{%
	postulates} to be satisfied by an area real valued function on polygons
which is supposed to be given. Then Moise [1964, pp.     165-7] proves this
postulate from its special case of a square of unit side, yet making use of
nontrivial properties of the reals. One can hardly say that any of these
alternatives would help Moschovakis.

On the contrary the Eudoxean approach proceeds as follows. Practically
speaking, the measure of the area of a polygon $P$ in terms of a (unit
square) tile $T$ is obtained as the result of comparison of $P$ with $T$. In
laying off $T$ on $P$, or vice versa, any of them may be reoriented, broken
and regrouped, or ... [cf. Lebesgue [1966] pp.     19, 42]. So, in fact what is
dealt with are equivalences of $T$ or $P$. The equivalence relation is
inspired by cns $1-4$ as shown in (iv) below [cf. Borceux [2014] p.    61].

Lebesgue [1966, p.    65] says ``one must have the concept of area before
calculating [measuring] areas''. This is achieved via kind (iv) which,
regarding its definition, embodies the concept of area of polygons. This
what makes Borceux [2014, p.    61] say ``In a sense, the Greek approach [to
area] is more ``intrinsic'' than ours, because it does not depend on the
choice  [beforehand] of a unit to perform the measure, and of course, different choices of
unit yield different measures of the same area.'' An important aspect of the
ingenuity of the Eudoxean approach is that it provides a universal definition
of (direct) measure. Once a kind is defined, the measure of one of its elements is the ratio of this element  to a specified
(unit) element. Based on this, the ratio  form of the well known formula for
the area of a rectangle (hence that of a triangle) may be deduced (see  II.2 
below) without recourse to real numbers, limits, nor exhaustion (the
definition of kind (iv) makes use only of Euclid [1956]'s book I). Though
extremely deep, the Eudoxean approach is very natural and intuitively quite
acceptable.

Moreover, let $F$ be the set of all positively Riemann integrable
non-negative real valued functions on closed bounded real intervals, and let 
$F^{\gimel }$ be the set of the regions under their graphs. Along the lines
of (v) below, $F^{\gimel }$ may be shown to be the universe of a kind $%
\mathbf{F}^{\gimel }$. With regard to $\mathbf{F}^{\gimel }$ it may be shown
that the integral of any of these functions is the positive real number
corresponding (see  II.2  below) to the ratio  of the corresponding region to
the unit square. This would provide Moschovakis with the needed
justification. Furthermore this shows why the area function may be real
valued.

In view of the \textit{current}, not only the historic, strong, probably
matchless, conceptual explanatory power of the Eudoxean approach, the
following thesis may be proposed: for each of the kinds considered below,
the corresponding Edoxean direct measure (as here explicated) is
conceptually sound; any other measure on the same kind or a subset thereof
is to be assessed with reference to it. Those who may have reservations
about adopting this thesis are invited to accept it as a working hypothesis.
In the extreme case of complete rejection, this article may still be of
significance as an attempt to precisely explicate a notion which had
prevailed for more than two millennia, and to discuss its role in the
evolution and applications of mathematics (see below). At any rate, further discussion
towards deeper analysis is always welcome.

Following is a discussion of the five aforementioned kinds; denote them by  $"\mathbf{K_1}", \cdots,"\mathbf{K_5}$", respectively.
 
\begin{enumerate}
	\item[(i)] For slss, $K$ is the set of all congruence classes of slss. The
	definitions of $+,<$ and $\dot{-}$ are obvious. The Archimedean
	principle is assumed.
	Notice that for the Archimedean principle to hold in each of cases (ii) -
	(v) discussed below, it is sufficient that it holds for the slss.
	\item[(ii)] For arcs of congruent circles, $K$ is the set of their congruence
	classes. For the definition of $+$, multiples of complete circumferences
	also should be allowed.
	\item[(iii)] The treatment of angles parallels that of arcs of congruent circles.
	\item[(iv)] The case of polygons is more involved. Hilbert [1950, p.    37]
	distinguishes between two equivalence relations on polygons: of \textit{%
		equal area} and of \textit{equal content}, which will be denoted here by ``$\rho
	_{0}$'' and ``$\rho _{1}$'' respectively [cf. Euclid [1956] vol. 1, p.    328]. 
	
	For a pair of polygons $P,P^{\prime },P\rho _{0}P^{\prime }$ iff they are
	piecewise congruent; while $P\rho _{1}P^{\prime }$ iff there are two
	polygons $P_{1}$ and $P_{1}^{\prime }$ such that:
	\begin{description}
		\item  (a) $P$ and $P_{1}$ are essentially-disjoint. Same for $P^{\prime }$
		and $P_{1}^{\prime }$.
		\item  (b) $P_{1}\rho _{0}P_{1}^{\prime }$ and $(P\cup P_{1})\rho
		_{0}(P^{\prime }\cup P_{1}^{\prime })$.
	\end{description}
	
	The definitions of $\rho _{0}$ and $\rho _{1}$are guided by Euclid's common
	notions [see above]. The definition of $\rho _{0}$ is guided by cn 4 "Things which coincide [congruent] with one another are equal to one another" and cn
	2 ``If equals are added to equals, the wholes are equal''. The definition of 
	$\rho _{1}$ is guided, in addition, by cn 3 ``If equals are subtracted from
	equals, the remainders are equal''. Obviously $\rho _{0}\subseteq \rho _{1}$.
	
	If $K$ is taken to be the $\rho _{0}$-equivalence classes, the definition of
	$+$, hence of $<$ and $\dot{-}$, will present no problem,
	but the trichotomy of $<$ will be problematic. Replacing $\rho _{0}$ by $%
	\rho _{1}$ will solve the problem. For, from proposition I.45 of Euclid and
	its proof [Euclid [1956] vol. 1, pp.     345-6] one may easily deduce: 
	
	{COROLLARY E}
	For every sls $l$ and every polygon $P$, there is a rectangle 
		$R(l,P)$ which is $\rho _{1}$-related to $P$ and which has a side congruent
		to $l$.

	The uniqueness (up to congruence) of $R(l,P)$ is guaranteed by cn
	5 "The whole is greater than the part" [cf. the proof of proposition I.39 of Euclid [1956, vol. 1, p.    336]]. From
	this the trichotomy of $<$ readily follows by the corresponding property of
	kind (i). So $K$ may be taken to be the $\rho _{1}$-equivalence classes.
	
	Why not to go beyond $\rho _{1}$? To be sure, guided by Euclid's cns, a
	monotonic increasing sequence $<\rho _{n}>_{n\in \Bbb{N}}$ of equivalence
	relations may be recursively defined as follows. $\rho _{0}$ is already
	defined. For every $n\in \Bbb{N}$, $\rho _{n+1}$ is to be obtained from $%
	\rho _{n}$ in a way similar to that by which $\rho _{1}$ is obtained from $%
	\rho _{0}$ above. The union of this sequence, $\sigma $ say, also is an
	equivalence relation.
	
	As a matter of fact, the proof of proposition I.45 of Euclid [1956, vol. 1,
	pp.     345-6] did not go beyond $\rho _{1}$. Moreover $\sigma =\rho _{1}$.
	Following is a proof based on cns 2-5 and the terminology thereof.
	
	By cns 4, 2, for every pair of polygons $P,P^{\prime }$, if $P\rho _{0}P^{\prime }$,
	then $P$ and $P^{\prime }$ are equal.
	Let $n\in \Bbb{N}$ and assume that for every pair of polygons $P,P^{\prime
	}, $ if $P\rho _{n}P^{\prime }$, then $P$ and $P^{\prime }$ are equal. Then
	by cn 3, for every pair of polygons $P,P^{\prime }$, if $P\rho
	_{n+1}P^{\prime } $, then $P,P^{\prime }$ are equal.
	
	By mathematical induction and the definition of $\sigma $, it follows that
	for every pair of polygons $P,P^{\prime }$, if $P\sigma P^{\prime }$, then $%
	P $ and $P^{\prime }$, are equal.
	
	Let $P,P^{\prime }$, be two polygons such that $P\sigma P^{\prime }$, and
	let $l$ be a sls. By corollary E, $R(l,P)\rho _{1}P\sigma P^{\prime }\rho
	_{1}$ $R(l,P^{\prime })$, hence $R(l,P)\sigma R(l,P^{\prime })$,
	consequently $R(l,P)$ and $R(l,P^{\prime })$ are equal. So, by cn 5 none of
	them may be a part of the other. From this it follows that they are
	congruent, then $R(l,p)\rho _{0}R(l,P^{\prime })$, a fortiori $R(l,P)\rho
	_{1}R(l,P^{\prime })$. Summing up: 
	
	$P\rho _{1}R(l,P)\rho _{1}R(l,P^{\prime
	})\rho _{1}P^{\prime }$. By the transitivity of $\rho _{1}$, $P\rho
	_{1}P^{\prime }$, which completes the proof.
	
	Making use of the Archimedean principle (for slss), Lebesgue [1966, p.    58]
	strengthened corollary E showing that $\rho _{1}$ may be replaced by $\rho
	_{0}$. So, by a simplification of the above argument, $\sigma =\rho _{0}$,
	and $K$ may be taken to be the $\rho _{0}$-equivalence classes.
	
	It is worth noting that no recourse to the real numbers is needed here. For
	other views see Hilbert [1950, pp.     37-45] and Lebesgue [1966, pp.     42-62]. We
	shall come back to this at the end of this section.
	\item[(v)] Regions bounded by slss and circular arcs: the set $C$ of all of these
	regions may be defined as follows. First, recursively define $<A_{n}>_{n\in 
		\Bbb{N}}$ by:
	\begin{description}
		\item  $A_{o}$= the set of all triangles $\cup $ the set of all circles.
	\end{description}
	 
	  $A_{n+1}=A_{n}\cup \{X$: there are $Y,Z\in A_{n}$ such that the
	set theoretic difference
	$Y-Z$ has a non-empty interior and $X$ is the (topological) closure of
	$Y-Z\}$.
	
	Put:
	\begin{description}
		\item  $B= \displaystyle \bigcup_{n \in \Bbb{N}} A_{n}$
		\item  $C$ = the set of all finite non-empty essentially-disjoint unions of
		elements of $B$.
	\end{description}

	Every element of $B$ is regular closed (i.e. is the closure of its
	interior), hence perfect (i.e. equals the set of all its limit points). Same
	applies to every element of $C$. The intersection of two elements of $B$ may
	not belong to $B$, but if the interior of this intersection is non-empty,
	then its closure belongs to $B$.
	
	Let $\eta $ be the binary relation defined on $C$ by: $X\eta Y$ iff for
	every inner polygon $P_{X}^{i}$ of $X$ and every outer polygon $P_{Y}^{o}$
	of $Y$: 
	
	\[
	P_{X}^{i}/\rho _{1}\leq P_{Y}^{o}/\rho _{1} 
	\]
	where $<$  is the ordering relation defined in (iv) above.
	
	$\eta $ is reflexive and, by the principle of exhaustion (which is due to
	Eudoxus and is equivalent to the Archimedean principle [Euclid [1956] vol.
	3, pp.     14-6, 365-8, 374-7 and vol. 1, p.    234]) for $C$, $\eta$ is transitive.
	I.e. $\eta $ is a pre-ordering, hence it defines an equivalence relation on $%
	C$, to be denoted by ``$\zeta $'', by: $X\zeta Y$ iff $X\eta Y$ and $Y\eta X$. $K$
	may be taken as the set of all $\zeta $-equivalence classes.
\end{enumerate}

Notice that the Archimedean principles for $C$ and slss are equivalent.
Assuming any of them, $\rho _{1}$ may be replaced by $\rho _{0}$ and the
mapping $P/\rho _{0}\longmapsto P/\zeta $, where $P$ is a polygon, is an
embedding of the kind defined in (iv) into the kind here defined.
%\end{itemize}

For each kind, direct measure [see  I.5  above] may be defined. The direct
measure of an element of any of the kinds (i)-(v) may be considered to be at
the same time the direct measure of each element of this element. Following
what may be called the Euclidean tradition, we may confuse equivalence
classes with elements thereof; the intention will be clear from the context.

The advantage of Eudoxus' treatment over that of Lebesgue is that it provides a unified approach to direct measure of all kinds, creating real numbers, while Lebesgue's treatment [Lebesgue [1966] pp.     19-67] is diversified, it is (with few exceptions) indirect, pre-assuming real numbers. 

Going back to Hilbert, it is worth noting that Hilbert [1950, p.    37] takes
into consideration not only $\rho _{0}$, but also $\rho _{1}$ ``to establish
Euclid's [Eudoxus'] theory of areas by means of the axioms already
mentioned; that is to say, {\it{for the plane geometry, and that is independent of
		the axiom of Archimedes.}}'' Nevertheless, this promise is not completely
fulfilled. Hilbert [1950] does not establish a theory of areas for kind (v),
it does [pp.     37-45] only for kind (iv) taking the well-known formula for the area of a triangle as a definition [p.    41] after proving that it does not depend on the choice of the base, without any further justification.

Likewise, without assuming the axiom of Archimedes, Hilbert [1950, pp.    
23-36] establishes a theory of proportion; it deals only with kind (i), in
particular, it does not deal with kind (iii): angles.

It is noteworthy that Tarski [1959, p.    21] adopted this theory of
proportion. Unlike that of Hilbert, Tarski's system of geometry is one
sorted which makes it better suited for metamathematical investigations.

As a matter of fact, Hilbert [1950, pp.     23-45]'s theories of proportion and areas are obtained
from the corresponding theories of Eudoxus by some sort of reverse
engineering (without mentioning Eudoxus). Computationally, the corresponding
theories are equivalent (assuming the axiom of Archimedes); conceptually,
they are not. In particular, from the measure-theoretic viewpoint, those of
Hilbert are entirely conceptually unsound as they disregard the concept of
direct measure. Moreover, unlike Euclid [1956], Hilbert [1950] and Tarski
[1959] do not touch the problem of measuring angles.

The serious defect of the aforementioned $20^{th}$ and $21^{st}$ centuries  approaches to trigonometry  [see
 I.1  above] is that they ignore altogether that there is a direct measure for
angles [see  II.4, II.5  below]. Hardy and Morrey [cf.  I.1  above] ignore,
moreover, that there is a direct measure for the plane regions they are
concerned with. The axiomatic and the purely analytic approaches to defining
the measure of angles or the trigonometric functions do not solve the
problem, they just transfer it [cf. the quotation from Moschovakis in  I.4 
above].

\vskip0.5cm
\textbf{II.2 Ratios and Real Numbers.} 
What are the real numbers? Over the past century (and, possibly, few more decades) a typical reply would be one of the following three options: equivalence classes of Weierstrass, nested closed intervals of rational numbers, equivalence classes of Cantor's fundamental sequences of rational numbers or Dedekind cuts of rational numbers [cf. Hobson [1957, p.    22]. This diversity is -in some sense- only apparent.

  All options define isomorphic models of one and the some theory:   Archimedean ordered fields with continuity. This theory is categorical, in the sense that \textit{all} of its models are isomorphic. 

  Where does this theory come from? A reply may be found in point (6) below.

  All of the three options presuppose the rational numbers, which historically originated from the measurement of (commensurable) magnitudes [cf. Hobson [1957, p.    14].

  How did Dedekind proceed stating from the rationals?
  His reply [Dedekind [1901, pp.     1-14] may be summarized as follows:
  
  \begin{itemize}
  	\item[1-]  Making use of the ratioal relations [pp.     3,4] he injectively maps the rational numbers into the set of points of a straight line, preserving order. After noticing that in the straight line there are infinitely many points which correspond to no rational number, he says that, if now, as is our desire, we try to follow up arithmetically all phenomena in the straight line it is necessary that the rationales be improved by newly \textit{created }numbers ``such that the domain of numbers shall gain the same [{\dots}] \textit{continuity}, as the straight line.''[p.    4].
  	\item[2-]  	However, ``this surely is no sufficient ground for introducing these foreign [geometric] notions into arithmetic, the science of numbers. Just as negative and fractional rational numbers are formed by a new \textit{creation} [emphasis added], and the laws of operating with  these numbers must and can be reduced to the laws of operating with positive integers, so we must endeavor completely to define irrational numbers by means of the rational numbers alone. The question only remains how to do this.'' [p.    5].
  	\item[3-]  	In fact, Dedekind considers arithmetic to be ``that part of logic which deals with the theory of numbers. In speaking of arithmetic (algebra, analysis) as a part of logic I mean to imply that I consider the number-concept entirely independent of the notions or intuitions of space and time, that I consider it an immediate result from the laws of thought.'' [p.    14].
  	
  	He adds, ``It is \textit{only} [emphasis added] through the purely logical process of building up the science of numbers and by thus acquiring the continuous number-domain that we are prepared accurately to investigate our notions of space and time by bringing them into relation with this number-domain created in our mind.'' [p.    14].
  	\item[4-] Before that, ``the way in which the irrational numbers are usually introduced is based directly upon the conception of extensive magnitudes -which itself is nowhere carefully defined- and explains number as the result of measuring such a magnitude by another of the same kind.'' [pp.     4-5].
  	\item[5-] Dedekind says that he finds the \textit{essence} of the continuity of the straight line in the following principle: Every cut is produced by exactly one point. [p.    5]. Moreover, he says that this property of the line, by which ``the secret of continuity is to be revealed.'' [p.    5], was taken for granted at his time [p.    5]. As he finds no proof of it, he takes it as an axiom ``by which we find continuity in the line.'' [p.    5]. call this axiom ``the axiom of continuity''.
  	\item[6-]  	So, we may safely say that at Dedekind's time it was implicitly taken for granted that the geometric field of real numbers is an Archimedean ordered field with continuity. Call it the ``Eudoxus-Descartes field of real numbers" or, for short, the ``E-D field''. 
  	
  	The contribution of Dedekind is that he explicated the notion of continuity. Thus, he practically informally axiomatized the theory of the E-D field; it is the (categorical) theory of Archimedean ordered fields with continuity, of which the E-D field is a geometric model. But Dedekind did \textit{not} develop a new theory.
	\item[7-] Via his cuts in the rational numbers, Dedekind \textit{creates} his irrational numbers, defining order on his real numbers (rationales and irrationals) and operations on them in such a way that practically makes what may be called the Dedekind field of real numbers isomorphic to the E-D field, hence it is a (n arithmetical) model of the theory of Archimedean ordered fields with continuity.
   \end{itemize}

  Commenting on Dedekind's \textit{creation }of irrational numbers, Frege [2016, vol. II, p.    141] says: ``This creating is the heart of the matter. [{\dots}]. Here the question is whether creating is possible at all; whether, if it is possible, it is so without constraint; or whether creation laws have to be obeyed while creating. 
  In the last case, before one could carry out an act of creation, one would first have to prove that, according to these laws, the creation is justified. These examinations are entirely missing 
  here, and thus the main point is missing; what is missing 
    is that on which the cogency of proofs conducted using irrational numbers depends.'' Another criticism of Dedekind's and other approaches of his contemporaries is, according to Frege [2016, vol. II, p.    157], that in these approaches: ``either measurement does not feature at all, or that it features without any internal connection grounded in the nature of the number itself, but is merely tacked on externally, from which it follows that we would, strictly speaking, have to state specifically for each kind of magnitude how it should be measured, and how a number is thereby obtained. 
    Any general criteria for where the numbers can be used as measuring numbers and what shape their \textit{application } [emphasis added] will then take, are here entirely lacking.'' As for Frege [2016, vol. II, p.    100], ``it is applicability alone which elevates arithmetic above a game to the rank of a science. Applicability thus necessarily belongs to it. Is it appropriate, then, to exclude from arithmetic what it needs to be a science?'' In conformity with this Frege [2016,vol. II, p.    85] sees that ``\textit{real number} is the same as magnitude-ratio.''

  Meanwhile  [2016, vol. II, p.    97] supported an important goal of the approaches of Dedekind and his contemporaries:  "Now, the real goal that all these theories of the irrationals are meant to serve is of course to display arithmetic as free from all foreign additions, including geometrical ones; to ground it on logic alone. This goal is certainly a proper one; [{\dots}].''

  After the crisis of Russell's paradox, late papers of Frege ``show him attempting a geometrical foundation for arithmetic.'' [Simons [1992, p.    137].

  Can we now answer the question we started with: what are the real numbers? Partially yes. 
  They are the elements of a model of the categorical theory of Archimedean ordered fields with continuity. Which one? Your choice!

  The choice in the rest of this article is the E-D field, not only because of its historical importance, but also because it is needed to provide a conceptually sound and logically correct foundations for trigonometry and analytic Euclidean geometry.
\vskip0.5cm
\textbf{II.2.1. The Eudoxus-Descartes Field of Real Numbers (The E-D Field).} 
The definitions of the operations on this field given by Descartes [1954] are problematic. He says [p.    5], ``For example, let $AB$ be taken as unity, and let it be required to multiply $BD$ by $BC$. I have only to join the points $A$ and $C$, and draw $DE$ parallel to $CA$; then $BE$ [which is to $BD$ as $BC$ to $AB$ [cf. p.    2]] is the product of $BD$ and $BC$.'' $BE$ depends not only on $BD$ and $BC$, but also on $AB$, the unity, ``which can in general be chosen arbitrarily, {\dots}." [p.    2].

  Is the product arbitrary? Is it relative (to the unity), not absolute? Is Descartes [1954] actually seeking the product of $BD$ and $BC$, or -to make the product absolute- he is seeking the product of the ratio of $BD$ to the chosen unity $AB$ and the ratio of $BC$ to the same unity? What are these ratios? Recall that De Morgan [1836, p.    62] noticed that the treatment of ratios by Euclid [1956] is problematic.

  It seems that the basic problem here is the definition of ratios, showing that they are endowed with the properties we intuitively ascribe to them. So, let us begin with it.

\vskip0.5cm
\textbf{II.2.1.1. Ratios.} 
The definition of ratios here parallels that of defining cardinal numbers in set theory. 
First, ``equinumerous'' is defined, then it is proved to be an equivalence relation. The most obvious continuation in to define the cardinal numbers as the equivalence classes of this equivalence relation. The problem faced is that these equivalence classes are not sets, but proper classes. To get around this problem,  a subset or an element is systematically chosen from each class, usually the smallest ordinal belonging to it. Via the well ordering principle, it is shown that this choice is possible and that the cardinal numbers so defined are endowed with all of the properties they are required to enjoy. 

  Likewise, concerning ratios the starting step, after making the notion of a kind (of magnitude) precise (see II.1 above), is to define "ratioal equivalence" which is given by:

DEFINITION V.5 of Euclid [1956, vol. 2, p.    114]
Magnitudes are said to 
\textbf{be in the same ratio  } [ratioally equivalent or proportional], the first to the
second and the third to the fourth, when, if any [positive integral]
equimultiples whatever be taken of the first and the third, and any
[positive integral] equimultiples whatever of the second and the fourth, the
former equimultiples alike exceed, are alike equal to, or alike fall short
of, the latter equimultiples respectively taken in corresponding order.

  Denoting this relation by ``$=_E$'', the above definition may be understood as follows. For  $x_1$,\textit{ }$x_2$ of the same kind, and $x_3$, $x_4$ of the same kind, $<x_1$,\textit{ }$x_2>=_E<x_3$,\textit{ }$x_4>$ iff for every pair of positive integers $m,n;$ $%
  mx_{1}>$ (respectively = or $<$) $nx_{2}$ iff $mx_{3}>$
  (respectively = or $<$) $nx_{4}$.

  That $=_E$ is reflexive and symmetric is obvious, proposition V. 11 of Euclid [1956, Vol. 2, p.    158] shows that it is transitive; so it is an equivalence relation. 

  Recall that we follow the old tradition which goes back to the ancient Greeks [II. 1 above] and which is most emphasized by Frege [2016,vol. II, p.    156]: ``Rather, we have exactly the same magnitude-ratio occurring in the case of periods of times, masses, light intensities, etc., as in the case of line segments.'' This is to be understood in the sense that there are $x_1$,\textit{ }$x_2$ of the same kind (which is different from kind ($i$) above) and two elements $y_1$,\textit{ }$y_2$ of kind ($i$) such that $<x_1$,\textit{ }$x_2>=_E$$\mathrm{<}$$y_1$,\textit{ }$y_2$$\mathrm{>}$. 
  In general, for "$<x_1,x_2>=_E<y_1,y_2>$" to hold, the $x$'s and the $y$'s do not have to be of the same kind.
   
  On the other hand there could be $x_1$,\textit{ }$x_2$ of the same kind such that for every $y_1$,\textit{ }$y_2$ of kind ($i$), $<x_1$,\textit{ }$x_2>{\neq }_E$$\mathrm{<}$$y_1$,\textit{ }$y_2$$\mathrm{>}$.

  So, taking all kinds into consideration, equivalence classes of $=_E$ may be proper classes.

  A good way out is to find a kind $\mathbf{K}$ which is universal in the sense that for each kind $\mathbf{K'}$, for each $<x',y'>\in \left|\mathbf{K'}\right|$ there is  $<x,y>\in \left|\mathbf{K}\right|$ such that $<x',y'>=_E<x,y>$. In this case the ratios may be defined as the equivalence classes of the equivalence relation $\left(\left(\left|\mathbf{K}\right|\times\left|\mathbf{K}\right|\right)\right)\times\left(\left(\left|\mathbf{K}\right|\times\left|\mathbf{K}\right|\right)\right)\cap =_E$, as these equivalence classes are sets.

  The history, and the current practice, of science nominate the kind $\mathbf{K_1}$ to be universal. For this to be the case it is sufficient that continuity holds for straight lines. To see this we go through the following steps. 

  (1) Following Tarski [1959], choose two different points $z$ (for zero) and $u$ (for unit). These two points determine a straight line $L(z, u)$. Via the betweenness relation, Tarski [1959, pp.     21-2] defines the familiar total ordering on the points of $L(z,u)$; denote  the corresponding strict total ordering (sto) by ``$<_{z,u}$''. For example,  $z<_{z,u}u$.

(2) On the other hand, toward totally ordering the equivalence classes of $=_E$, consider:

   DEFINITION. V.7 of Euclid [1956, vol. 2, p.    114]. When, of the equimultiples,
   the multiple of the first magnitude $[x_1]$ exceeds the multiple of the second $[x_2]$, but
   the multiple of the third $[x_3]$ does not exceed the multiple of the fourth $[x_4]$, then
   the first is said to \textbf{have a greater ratio} to the second than the
   third has to the fourth [$<x_1,x_2>$ ratioally dominates $<x_3,x_4>$].

  Denoting this relation by ``$>_E$'' (and its converse by ``$<_E$''), the above definition may be understood as follows.

  For $x_1$,\textit{ }$x_2$ of the same kind and $x_3$,\textit{ }$x_4$ of the same kind, $<x_1$,\textit{ }$x_2>>_E<x_3$,\textit{ }$x_4>$ (equivalently, $<x_3,\ x_4>\mathrm{\ }<_E<x_1$,\textit{ }$x_2>$) if the above condition, i.e., there are positive natural numbers m, n such that $mx_1>nx_2$ and $mx_3\le nx_4$, is satisfied. In this case we say that $<x_1$,\textit{ }$x_2>$ ratioally dominates $<x_3$,\textit{ }$x_4>$, and $<x_3$,\textit{ }$x_4>$ is ratioally dominated by $<x_1$,\textit{ }$x_2>$. 

  It is easy to see that $<_E$ is a strict partial ordering, but it is not total, though it is close to being so. For it is not hard to see that $<x_1$,\textit{ }$x_2>{\nless }_E<x_3$,\textit{ }$x_4>$ and $<x_3$,\textit{ }$x_4>{\nless }_E<x_1$,\textit{ }$x_2>$ iff $<x_1$,\textit{ }$x_2>=_E<x_3$,\textit{ }$x_4>$. So, by invoking proposition V.13 of Euclid [1956, vol. 2, pp.     160-1] and a variation thereof, it may be easily seen that $<_E$ induces a strict total ordering (sto) on the equivalence classes of $=_E$. The induced relation will be denoted by ``$<'_E$''.

  Alternatively, Eudoxus could have defined a binary relation, to be denoted by "$\circleit{$\scriptscriptstyle <$}$", as follows:

  For $x_1$,\textit{ }$x_2$ of the same kind and $x_3$,\textit{ }$x_4$ of the same kind,  $<x_1$,\textit{ }$x_2>$ $\circleit{$\scriptscriptstyle <$}$ $<x_3$,\textit{ }$x_4>$ if there are positive natural numbers $m,\ n$ such that $mx_1<nx_2$ and $mx_3\ge nx_4$.

  Like above, we may show that $\circleit{$\scriptscriptstyle <$}$ is a strict partial ordering, moreover ,  $<x_1$,\textit{ }$x_2>$ $\circleit{$\scriptscriptstyle \not<$}$ $ <x_3$,\textit{ }$x_4>$ and $<x_3$,\textit{ }$x_4>$ $\circleit{$\scriptscriptstyle \not<$}$ $<x_1$,\textit{ }$x_2>$ iff $<x_1$,\textit{ }$x_2>=_E<x_3$,\textit{ }$x_4>$.

  In fact, $\circleit{$\scriptscriptstyle <$}$ $=<_E$. To see this, show that , ($<x_1$,\textit{ }$x_2><_E<x_3$,\textit{ }$x_4>\mathrm{,\ \ }<x_1,\ x_2>$ $\circleit{$\scriptscriptstyle >$}$ $<x_3,\ x_4>$) leads to contradiction. So $<x_1$,\textit{ }$x_2><_E<x_3$,\textit{ }$x_4>$ implies $<x_1$,\textit{ }$x_2>$ $\circleit{$\scriptscriptstyle \not>$}$ $<x_3$,\textit{ }$x_4>,$\textit{ } which implies ($<x_1$,\textit{ }$x_2> $ $\circleit{$\scriptscriptstyle <$}$ $<x_3$,\textit{ }$x_4>$ or $<x_1$,\textit{ }$x_2>=_E<x_3$,\textit{ }$x_4>$), but $<x_1$,\textit{ }$x_2><_E<x_3$,\textit{ }$x_4>$ implies $<x_1$,\textit{ }$x_2>{\neq }_E<x_3$,\textit{ }$x_4>$, hance $<x_1$,\textit{ }$x_2><_E<x_3$,\textit{ }$x_4>$ implies $<x_1$,\textit{ }$x_2>  $ $\circleit{$\scriptscriptstyle <$}$ $<x_3$,\textit{ }$x_4>$. 

  The other direction may be proved similarly. 
  
  As above, $\circleit{$\scriptscriptstyle <$}$ induces a strict total ordering on the equivalence classes of $=_E$, to be denoted by "$\circleit{$\scriptscriptstyle <$}'$". So $\circleit{$\scriptscriptstyle <$}'$ $=<'_E$.

(3) Recall that $K_1=\left|\mathbf{K_1}\right|$ and $L(z, u)$ is the straight line determined by points $z, u$. Put:

  (a) $=_{K_1}=(K_1\times K_1)\times (K_1\times K_1)\cap =_E$, then $=_{K_1}$ is an equivalence relation on $(K_1\times K_1)$,

  (b)$<'_{K_1}=((K_1\times K_1)/=_E\times{\left(K_1\times K_1\right)/=}_E)\cap <'_E$, then $<'_{K_1}$ is a sto on $(K_1\times K_1)/=_E$, and via the bijection $(<x,y>/=_{K_1}\longmapsto <x,y>/=_E)$ from $(K_1\times K_1)/=_{K_1}$ onto $(K_1\times K_1)/=_E$ ,\textit{  }$<'_{K_1}$ induces a sto on $(K_1\times K_1)/=_{K_1}$, to be denoted by ``$<''_{K_1}$''.

  Also, put:

  (c)$Ray\left(z,u\right)=\left\{x\in L\left(z,u\right):z<_{z,u}x\right\}$,

  (d)$Rs\left(z,u\right)=\left\{\overline{zx}:x\in Ray(z,u)\right\}$, where $\overline{zx}$ is the sls with end points $z,x$,

  (e)$Rs\left(z,u\right)/cong=\left\{\overline{zx}/cong:\overline{zx}\in Rs(z,u)\right\}$, where "\textit{cong}" is the well known geometric congruence relation. But, 

  (f)$Rs$$\left(z,u\right)/cong=K_1,$ so $ \ <''_{K_1}$ is a sto on ($Rs\left(z,u\right)/cong\ \times\ Rs\left(z,u\right)/cong)/=_{K_1}(=(K_1\times K_1)/=_{K_1})$

  (g)By the existence of the 4${}^{th}$ proportional for slss (proposition VI.12 of Euclid [1956, vol. 2, p.    215]):

  ($Rs\left(z,u\right)/cong\ \times \ \left\{\overline{zu}/cong\right\})/=_{K_1}=\mathrm{(}Rs\left(z,u\right)/cong\ \times  \ \mathrm{(}Rs\left(z,u\right)/cong)/=_{K_1}$, hance $<''_{K_1}$ is a sto on ($Rs\left(z,u\right)/cong\ \times \ \left\{\overline{zu}/cong\right\})/=_{K_1}$ 

  (h)By proposition V.9 of Euclid [1956, vol.2, pp.    153-4] there is a bijection $(\overline{zx}/cong\longmapsto <\overline{zx}/cong,\overline{zu}/cong>/=_{K_1})$ from $Rs(z,u)/cong$ onto $(Rs(z,u)/cong\ \ \times \ \left\{\overline{zu}/cong\right\})/=_{K_1}$. Also, there is an obvious bijection from $Rs\left(z,u\right)\ onto\ Rs(z,u)/cong$ .
  Again, there is an obvious bijection from $Ray\left(z,u\right)$ onto $\ Rs(z,u)$ .
  So, $<''_{K_1}$ induces three sto relations (to be denoted, respectively, by $<''_{Rs},<'_{Rs}$ and $<_{Ray(z,u)}$) on $Rs(z,u)/cong\left(=K_1\right),\ Rs(z,u)$ and $Ray(z,u)$, which make the four order structures isomorphic.

  (i)\textit{ }$<''_{Rs}$ is the same sto given in the definition of the kind $\mathbf{K_1}$, $<'_{Rs}$ is the same sto defined by Hilbert [1950, p.    30], $<_{Ray(z,u)}$ is the restriction on $Ray(z,u)$ of  $<_{z,u}$ defined in (1) above.

(j) {LEMMA.}
 Let $s$ be a sls, $\mathbf{K}$ a kind, $v,v'\in K$ and $m,n\in \Bbb{N}^+$, then: 
	\[
	\begin{split}
	<ms,ns>&=_E<v,v'>\qquad\qquad\qquad \text{iff}\qquad mv'=nv\\&
	<_E\qquad\qquad\qquad\qquad\qquad \qquad\qquad\,\,\,\,<\\&
	>_E \qquad\qquad\qquad\qquad\qquad \qquad\qquad\,\,\,\, >
	\end{split}
	\]

\begin{proof}    
The proof of the first part is straightforward. For the second part, make use of the definition of  $<_E$. Show that one direction is straightforward. For the other direction show that the antecedent together with the negation of the consequent leads to contradiction.

  The proof of the third part in similar to that of second part, but making use of the definition of $\circleit{$\scriptscriptstyle <$}$ instead of the definition of $<_E$.
\end{proof}

(k) {LEMMA.}
 Let $\mathbf{K},\mathbf{K'}$ be kinds, and let $v,v_1\in K$, and $v',v'_1\in K'$. $<v,v_1><_E$ $<v',v'_1>$ iff  there are $i,j\in \Bbb{N}^+$ and an $s\in K_1$, such that:
\[<v,v_1>{<_E<is,js><}_E<v',v'_1>\] 
(Notice that for every $i,j\in \Bbb{N}^+$ and every $s,s'\in K_1$, $<is,js>=<is',js'>$)

\begin{proof}
   The if direction is obvious. To prove the only if direction, assume $<v,v_1><_E<v',v'_1>$. Distinguish between four cases: 

  Case1: for every $k,l\in \Bbb{N}^+,\ <v,v_1>\neq_E<ks,ls>\neq_E<v',v'_1>$.

  By the definition of $<_E$, there are $i,j\in \Bbb{N}^+$ such that $iv\le jv_1$ and ${iv}'>jv'_1$. Then, by the above lemma $ <v,v_1>{\le }_E<js,is><_E<v',v'_1>$. Hence $<v,v_1><_E<js,is><_E<v',v'_1>$.

  Case2: For same $k,l,k',l'\in \Bbb{N}^+$ and some sls $s,\ <v,v_1>=<ks,ls>$ and \\  $<v',v'_1>=<k's,l's>$. By the existence of the fourth proportion, there is a sls $s'$ such that:
\[<ks,ls>=_E<s',l's><_E<k's,l's>.\] 
By proposition V.10 of Euclid [1956, vol.2, pp.    155-6] $s'<k's$. By the Archimedean property and simple Euclidean geometry, there are $p,q\in \Bbb{N}^+$ and a sls $s''$ such that $s=qs''$ and $s'<ps''<k's$. By proposition V.8 of Euclid [1956, vol.2, pp.    149-52]:

  $<ks,ls>=_E<s',l's><_E<ps'',l's><_E<k's,l's>$, hence \\ $<ks,ls><_E<ps,ql's><_E<k's,l's>$.

  Case3: For some $k,l\in \Bbb{N}^+$ and some sls $s,<v',v'_1>=_E<ks,ls>$, while for every $k,l\in \Bbb{N}^+,\ <v,v_1>{\neq }_E<ks,ls>$.

  By the definition of $<_E$, there are $m,n\in \Bbb{N}^+$ such that:

  $mv\le nv_1$ and $mks>nls$. Hence:

  $nv_1\ge mv$ and $nls<mks$

  So, by the above lemma:
\[<v,v_1>\le_E<ns,ms><_E<ks,l,s>,\] 
Consequently:
\[<v,v_1><_E<ns,ms><_E<ks,ls>.\] 
Case4: For some $k,l\in \Bbb{N}^+$ and some sls $\ \ s,\ <v,v_1>=_E<ks,ls>$, while for every $k,l\in \Bbb{N}^+,\ <v',v'_1>{\neq }_E<ks,ls>$.

  To prove case4, follow the steps of proving case3, with making use of the definition of $\circleit{$\scriptscriptstyle  <$}$ instead of the definition of $<_E$.
\end{proof}

  (l) THE MAIN LEMMA. Let $\mathbf{K}$ be a kind, and $v,v'\in K$, and let $s,s'\in K_1$, then $\mathrm{<}$\textit{ }$v,v'$$\mathrm{>}$$=_E<s,s'>iff$
\[\sup\limits_{(K_1\times K_1)/=_{K_1}}\left\{<is,js>/=_{K_1}:i,j\in \Bbb{N}^+\ and\ <is,js><_E\mathrm{<}\ v,v'\mathrm{>}\right\}=<s,s'>/=_{K_1}\] 
  (in the sense that $\sup$ exists and equals the rhs).
 
\begin{proof}
   Assume $\mathrm{<}\ v,v'\mathrm{>}=_E<s,s'>$. Let $<is,js><_E\mathrm{<}\ v,v'>$, then $<is,js><_E<s,s'>$, then $<is,js>$/$=_{K_1}<''_{K_1}<s,s'>/=_{K_1}$, so $<s,s'>/=_{K_1}$ is an upper bound. Also let $s_1,\ s'_1\in K_1$, and  let $<s_1,\ s'_1>/=_K<''_{K_1}<s,s'>/=_{K_1}$, then ${<s}_1,\ s'_1><_E<s,s'>$. By the above lemma, there are $i,j\in \Bbb{N}^+$ such that ${<s}_1,\ s'_1><_E<is,js><_E<s,s'>=_E\mathrm{<}\ v,v'>$. Hence $<s_1,\ s'_1>/=_K$ is not an upper bound, so $ <s,s'>/=_{K_1}$ is the supremum. 

  On the other hand, let $<s,s'>/=_{K_1}$ be the supremum. Let $\mathrm{<}\ v,v'><_E<s,s'>$, then by the above lemma there are $k,l\in \Bbb{N}^+$ such that $\mathrm{<}\ v,v'><_E<ks,ls><_E<s,s'>$, hence $<ks,ls>/=_{K_1}$ is an upper bound and $<ks,ls>/=_{K_1}<''_{K_1}<s,s'>/=_{K_1}$. Which contradicts that $<s,s'>/=_{K_1}$ is the supremum. So $\mathrm{<}\ v,v'>\not<_E<s,s'>$.

  Also, let $<s,s'><_E\mathrm{<}\ v,v'>$, then by the above lemma there are $i,j\in N^+$ such that $<s,s'><_E<is,js><_E\mathrm{<}\ v,v'>$. Hence $<s,s'>/=_{K_1}$ is not an upper bound which contradicts that it is the supremum. So  $<s,s'>{\nless }_E\mathrm{<}\ v,v'>$. Therefore, $\mathrm{<}\ v,v'\mathrm{>}=_E<s,s'>$. 
\end{proof}

  (m) Most, but not all, of the propositions of Euclid [1956] made use of here apply to all kinds. In particular, the existence of the fourth proportional does not apply to the kind of positive natural numbers. So, what is proved here may apply (totally or partially) to kinds other than $K_1$, but we shall not further occupy ourselves by this matter here. 

  (n) The main lemma shows that the kind $K_1$ is universal (see above) if the supremum always exists. This is guaranteed if the ordering $<''_{K_1}$ on $(K_1\times K_1)/=_{K_1}$ satisfies the axiom of continuity.

  By part (h) above, this is guaranteed if the ordering $<_{Ray(z,u)}$ on $Ray(z,u)$ satisfies the axiom of continuity. 

  (o) As our aim is to precisely define ratios, toward precisely defining the E-D field, it is appropriate at this point to adopt the axiom of continuity for $<_{z,u}$, hence for $<_{Ray(z,u)}$, hence for $<''_{K_1}$.

  (p) THE MAIN DEFINITION / THEOREM.

  (i) A ratio is an element of $(K_1\times K_1)/=_{K_1}$(so, $\ (K_1\times K_1)/=_{K_1}$ is the set of all ratios). 

  (ii) Let $\mathbf{K}$ be a kind, and let $v,w\ \epsilon K$, then Ratio ($v,w$), read:
  the ratio of \textit{ }$<v,w>$, or the ratio of $v$ to $w$, is defined as follows:

  Ratio
  $(v,w)=\sup\limits_{(K_1\times K_1)/=_{K_1}}\{<is,js>/=_{K_1}:
  i,j\in \Bbb{N}^+\,\,\text{and}\,\,\, 
  <is,js><_E<v,w>\}$
   \\(the existence of the $sup$ is guaranteed by continuity).

  (iii) For a kind $\mathbf{K}$ and $v,w\ \in K$, there are $r,s\ \in K_1$ such that:
  $<v,w>=_E<r,s>$, so $\mathbf{K_1}$ is universal.

  (iv) Let $\mathbf{K},\ \mathbf{K}'$ be two kinds, and let $v,w\ \in K$ and $v',w'\ \in K'$, then:
\[
<v,w>=_E<v',w'>\,\,\text{iff}\,\,\,\, Ratio\,\,\left(v,w\right)=Ratio\,\,(v',w').
\]
   
\begin{proof}
   By continuity there are $r,s,r',s'\in K_1$ such that:

  Ratio $\left(v,w\right)=<r,s>/=_{K_1}$ and $Ratio\left(v',w'\right)=<r',s'>/=_{K_1}$. Then, by the main lemma:
\[<v,w>=_E<r,s>\ \mathrm{and\ }\ {<v',w'>=}_E<r',s'>.\] 
This proves (iii). Moreover, 

  $<v,w>=_E<v',w'>$ iff  ${<r,s>=}_E<r',s'>$ iff (by a, above)
\[<r,s>=_{K_1}<r',s'>\ \mathrm{iff\ }\ {<r,s>/}=_{{K_1}}\,\,\,\,=\,\,\,\,<r',s'>/=_{K_1}\] 
iff $Ratio\left(v,w\right)=Ratio(v',w')$.  This proves (iv).
\end{proof}

\vskip0.25cm
  \textbf{II.2.1.2. Operations.} Euclid [1956] has all needed ingredients to define two binary operations on $(K_1\times K_1)/=_{K_1}$ 
  (the set of all ratios). The first [cf. Euclid [1956] vol.2, p.   187 and Stein [1990] pp.    177, 181-3], to be called ``addition'' and to be denoted by ``$+$'', is defined as follows: For $<a,b>,\ <c,d>\in (K_1\times K_1)$, put:
\[<a,b>/=_{K_1}+<c,d>/=_{K_1}=<a+c',b>/=_{K_1}\] 
Where $c'(\in K_1)$ is such that $<c',b>=_{K_1}<c,d>$.

  This justifies the name ``addition'' and the symbol ``$+$'', though there is some abuse of notation. That + is well defined is guaranteed by propositions VI.12 of Euclid [1956, vol.2, p.   215] and V.24 of Euclid [1956, vol.2, p.   183]. It is not hard to see that + is associative and commutative and that $<_E$ is compatible with + .

  The second [cf. Euclid [1956] vol.2, p.   187 and Stein [1990] pp.    183, 185], to be called ``multiplication'' and to be denoted by ``$\times$'' (not to be confused with the Cartesian product, the intention will be clear from the context), is defined as follows:

  For $<a,b>,\ <c,d>\in (K_1\times K_1)$, put: 
\[<a,b>/=_{K_1}\times<c,d>/=_{K_1}\,\,\,\,=\,\,\,\,<a,d'>/=_{K_1}\] 
Where $d'$ is the unique element of $K_1$ such that :
\[<c,d>=_{K_1}<b,d'>\]
 
By propositions V22, 23 of Euclid [1956, vol.2, pp.    179-82] $\times$ is well defined, associative, commutative, distributive (over +) and $<_E$ is compatible with it. Notice that $\times$ is absolute, it does not depend on any explicit nor implicit parameters. Notice also that:

\[<a,b>/=_{K_1}\times<c,d>/=_{K_1}=<a',d>/=_{K_1}=<c,b'>/=_{K_1}<c',b>/=_{K_1},\] 
Where $a',\ b',\ c'$ are, respectively, the unique elements of $K_1$ such that 

\[<a,b>=_{K_1}<a',c>=_{K_1}<d,{\ b}'>,\ \mathrm{and}<c,d>=_{K_1}<c',a>.\] 

The above definition implies that:

\[
\begin{split}
Ratio(Rec'(a,b),Sq'(c))&=<a',c'>/=_{K_1}\times <b',c'>/=_{K_1}
\\&
=Ratio(a',c')\times Ration(b',c'),
\end{split}
\]
where $a,b,c$ are slss, and $a',b',c'$ are the elements of $K_1$ corresponding, respectively, to them; also, $Rec(a,b)$ is the rectangle whose sides are $a,b,Sq(c)$ is the square whose side is $c$, and ${Rec}'(a,b),\ {Sq}'(c)$ are, respectively, the elements of $K_4$ corresponding to them [cf. Euclid [1956] vol.2, prop.    VI.23, p.   247]. $Sq(c)$ may be looked upon as a unit square and $c$ the corresponding unit sls. 

  This justifies the choice of the name ``multiplication'' and the symbol ``$\times$''. 

\vskip0.5cm
  \textbf{II.2.1.3.Descartes, Hilbert and Tarski - Zero and the Negatives.} Tarski [1959, pp.    21.2] says that he is ordering, adding and multiplying points. For ordering there is no problem, only betweenness is made use of to order the points. For addition and multiplication matters are different. With some scrutiny we realize that Tarski is in fact adding slss. As all of his slss share one end point, namely $z$, each may be replaced by its other end point. For multiplication, Tarski implicitly refers the reader to a German edition of Hilbert [1950]. 

  Descartes [1954, pp.    2, 5] and Hilbert [1950, pp.    29-30] say that they are adding slss; there is no problem with this. Multiplication is problematic; they both say that what they are multiplying are slss, but only after arbitrarily choosing a third slss which remains fixed through the discussion; it may be called ``unit'' and denoted by ``$u$''. 

  Descartes [1954, pp.    2, 5] did not tell why did he define multiplication the way the did. What is the relationship between his definition [see the beginning of II.2.1 above] and the definition of multiplication of ratios given in II.2.1.2 above? To deal with these matters, let $a,b,p,u\ \mathrm{be,}$ respectively, the elements of $K_1$ corresponding to the two slss Descartes is multiplying, their product and the unit, then:
  
\[<p,b>=_{K_1}<a,u>,\]
then
\[<p,b>/=_{K_1}\times<b,u>/=_{K_1}=<a,u>/=_{K_1}\times<b,u>/=_{K_1},\] 
then 
\[<p,u>/=_{K_1}=<a,u>
/=_{K_1}\times<b,u>/=_{K_1},\] 
then 
\[Ratio<a,u>\times\ Ratio<b,u>=Ratio<p,u>.\] 

Therefore Descartes is, in fact, multiplying ratios. Form this and the definitions of $<''_{K_1}$ and + all the required properties of his multiplication (commutativity, {\dots}) follow at once. 

  Without adopting Eudoxu's theory of ratio and proportion, Hilbert [1950, p.   30] adopts Descartes [1954, pp.    2, 5] definition of multiplication of slss, without mentioning him. Why? Where did this definition come from? Is the product relative or absolute? What is to be done concerning the dependence on the unit $u$? He does not reply; In fact these questions were not raised to start with. 

  Hilbert [1950, pp.    25-9] states and proves Pascal's theorem without recourse to the Archimedean axiom, and he makes use of this theorem [pp.    31-2] to prove the required properties of multiplication. 

  In contrast to Descartes [1954, pp.    2, 5] who makes use of proportion to define multiplication, Hilbert [1950, p.   32] makes use of multiplication to define proportion, without discussing whether it is absolute or is relative to the unit $u$. 

  Beside implicitly referring the reader (for the definition of multiplication) to a German edition of Hilbert [1950], Tarski [1959, pp.    21-2] explicitly says ``The definition of $x.y$ is more involved; it refers to some points outside of $F$ and is essentially based upon the properties of parallel lines.'', ignoring that the definition depends on the unit $u$. 

  To complete the definition of the E-D field, we have to introduce zero and the negatives. To do so we follow -with some modifications- what Descartes implicitly did [cf. Descartes [1954] p.   196] and what Hilbert [1950, p.   34] and Tarski [1959, pp.    21-2] explicitly did.

  Go back to the straight line $L(z,u)$ defined in (1) of II.2.1.1 above. Consider all  slss with $z$ as one end point and the other end point is an arbitrary point of $L(z,u)$, possibly $z$ itself.

  So, $\overline{zz}$ will be considered here as a sls, which is very un-Greek. Also, we distinguish here between slss whose other end points are greater than $z$, which we call positive and denote by ``$P_{z,u}$'', and those whose other end points are smaller than $z$, which we call negative and denote by ``$N_{z,u}$''; again, this is un-Greek. 

  Put: $S_{z,u}=P_{z,u}\cup N_{z,u}\cup \left\{\overline{zz}\right\}$.

  Obviously, $Rs\left(z,u\right)=P_{z,u}$, hence: \\the set of all ratios  $=(K_1\times K_1)/=_{K_1}\,\,\,=\,\,\,(P_{z,u}/cong\ \times \left\{\overline{zu}/cong\right\})/=_{K_1}$

  Keeping $\overline{zu}$ (hence $\overline{zu}/cong$) in the background and noticing that for each element of $P_{z,u}/cong$ there corresponds exactly one element of $P_{z,u}$, the set of all ratios may be replaced by $P_{z,u}$. In fact to a ratio $r$ there corresponds the unique $\overline{zx}\in P_{z,u}$ such that :
  
\[r=Ratio(\overline{zx}/cong,\ \overline{zu}/cong)\] 
\[(=<\overline{zx}/cong,\ \overline{zu}/cong>/=_{K_1}).\] 

So we may consider + and $\times$ (defined on the set of all ratios)
  to be operations on $P_{z,u}$, and consider $<''_{K_1}$ to be a relation on $P_{z,u}$. 

  This relation together with the two operations may be extended in the obvious way to $S_{z,u}$ [cf. Tarski [1959] pp.    21-2] for order and addition, [cf. Frege [2016] vol.II, p.   160 for addition, and cf. Descartes [1954] pp.    2,5 for multiplication]. Thus $S_{z,u}$ becomes a continuously (hence Archimedean) ordered field. The unit elements for addition and multiplication are $\overline{zz}$ and $\overline{zu}$ respectively.

  Every positive $sls$ is congruent to exactly one negative $sls$ and vice versa, they are the additive inverses of each other. For division, hence the multiplicative inverse, Descartes [1954, pp.    2,5] may be consulted. Obviously  $P_{z,u}$ is the set of positive elements of $S_{z,u}$.

  For any two such fields, $S_{z,u}$ and $S_{z',u'}$, there is exactly one isomorphism from the former to the latter, namely the extension of the mapping which sends each positive element $\overline{zx}$ of $S_{z,u}$ to the unique positive element $\overline{z'x'}$ of $S_{z',u'}$ such that:
  
\[Ratio(\overline{zx}/cong,\overline{zu}/cong)=Ratio(\mathrm{\ }\overline{z'x'}/cong,\mathrm{\ }\overline{z'u'}/cong).\] 

Any of these fields may be regarded as the E-D field, or the (geometric) field of the real numbers. 

To make the E-D filed closer to the spirit of Eudoxus and Frege, the set of all ratios may be taken as the set of positive reals, and zero and the negatives be redefined in the obvious way. 
\vskip0.75cm
\textbf{II.2.1.4. Continuity. } Is adopting the axiom of continuity in general, and as an axiom of Euclidean geometry in particular, conceptually sound?

  The matter is problematic. An important consequence of this axiom is the formulation of the field of the real numbers which is problematic too. Its role in the scientific applications of mathematics is crucial, and at the same time we cannot claim that we deeply understand it; we do not know its cardinality so far.

  Can mathematics which is indispensable to scientific applications be developed by a weaker axiom than continuity? Possibly yes [cf. Feferman [1999] p.   109]. Is it practicable to expect that this yet to-be-developed mathematics will be as effective as the current mathematics? This is yet to be seen. Until answers to these questions are found, if the history of the subject is taken into consideration, I guess it is fair to regard the axiom of continuity as conceptually sound, at least for the lack of any better.

  Roughly speaking, adopting the axiom of continuity is equivalent to adopting all of the theoretically possible ratios. Frege, who is a strong advocate of regarding real numbers as ratios adopted continuity. He takes continuity as a condition of his definition of positive classes [Frege [2016] vol.II, p.   190]. On the positive classes are based the domains of magnitudes [cf. Frege [2016] vol.II, p.   170], and the real numbers are to be ratios of the elements of these domains [cf. Frege [2016] vol.II, pp.    155,157,160]. Frege's domains of magnitudes (or rather the positive classes themselves) are the counterparts of our kinds of magnitude. 

  What about the soundness of adopting continuity as an axiom of Euclidean geometry? Dedekind who formulated this axiom, said [Dedekind [1901] PP4-5] that he explicated it from Euclidean geometry. He adds [p.   5] that continuity was implicitly taken for granted by his contemporaries as an axiom of this geometry. So it is reasonable to adopt continuity as an axiom of Euclidean geometry, if is to be adopted at all. 

  Lebesgue [1966, pp.    20-1] also regards continuity to be an implicitly understood axiom of Euclidean geometry, though it is not explicitly stated.      

  A careful reading of Descartes [1954] shows that he probably had implicitly adopted continuity as an axiom of Euclidean geometry. 

  Hilbert [1950, p.   36] practically grants that continuity is an axiom of the ``ordinary analytic geometry of space''.

  Tarski [1959, p.   18] adopts the continuity axiom as a member of ``an axiom system which is known to provide an adequate basis for the whole of Euclidean geometry''. As this axiom is second order, not elementary, it was replaced by ``the infinite collection of all elementary continuity axioms'', in the axiomatization of elementary geometry adopted by Tarski [1959, pp.    18-20]. 

  Concerning continuity Stein [1990, p.   178] said that it ``is by no means necessary for the geometry contained in Euclid's \textit{Elements}- a fact already emphasized by Dedekind himself in the preface to the first edition of his monograph on natural numbers [Dedekind [1901] pp.    16-17].'' To vindicate his statement Stein [1990, pp.    178-9] refers to the smallest subfield of the field of the real numbers closed under taking the positive square root of the positive elements, and its relationship to the Elements. A detailed treatment of this (which is a variant of an argument advanced by Dedekind [1901, pp.    16-7]) may be found in Moise [1964, chap.    19, pp.    214-41].

  However in dealing with the question whether the Greek geometers would have accepted Dedekind's axiom of continuity once it had been stated? Stein [1990, p.   179] answers ``there are reasons -suggestive, and, I think, plausible, even if not conclusive- for'' believing that they would have accepted it. 

  In light of this discussion, I guess it is reasonable to regard geometric continuity as conceptually sound, until further notice (see above). Accordingly, henceforth the geometric axiom of continuity will here be adopted. 

\vskip0.5cm
\textbf{II.3. Angles.} 
Euclid's definition of a (rectilineal)
angle is given [Euclid [1956] vol. 1, p.    153] in two steps: 

DEFINITION I.8 of Euclid [1956, vol. 1, p.    153]. \textbf{A plane angle} is
the inclination to one another of two lines in a plane which meet one
another and do not lie in a straight line.

DEFINITION I.9 of Euclid [1956, vol. 1, p.    153]. And when the lines
containing the angle are straight, the angle is called \textbf{rectilineal}.

Modern writers try to avoid the obscurity of the term ``inclination'' by
defining the rectilineal angle (henceforth, angle) as the configuratioal relationn
formed by two straight lines [rays] having the end point in common [Hardy
[1967] p.    316], or as the locus consisting of two rays which have a common
end point and lie along different lines [Morrey [1962] p.    203].

These definitions share the following two characteristics:
\begin{description}
	\item  (i) They make use of infinite objects: (full) rays, which is alien to
	Greek geometry.
	\item  (ii) Different angles, according to these definitions, may be equal
	according to Euclidean geometry.
\end{description}

Noting that in Euclidean geometry ``equal'' is sometimes used where
``congruent'' or ``equivalent'' should have been used, it is reasonable to
retain the second characteristic. To avoid the first one, the following
definition is proposed:

  An angle is an order pair $<<a,b,c>,d>$ whose first component is an ordered triple of three pairwise distinct points $a,b,c,$ and whose second component is a point $d$ not belonging to $Ray(b,a)\cup Ray(b,c)$; $b$ is called the vertex of the angle.

  Intuitively, the ordered triple $<a,b,c>$ defines two angles, and $d$ determines exactly one of them in the obvious way [cf. Hilbert [1950] pp.    8-9 and De Morgan [1837] p.   15]. In case $Ray(b,a)$ and  $Ray(b,c)$ coincide, $d$ always determines one and the same angle, the angle which is not the zero angle, i.e., the perigon or the round angle. 

  The reason behind excluding the zero angle here, is to allow the congruence classes of angles to form a kind ($\mathbf{K_3}$, above). Subsequently, this angle will be taken into consideration. 
\vskip0.5cm
\textbf{II.4 Measure of Angles. } 
Though the notion of a ratio was problematic, as De Morgan [1836, p.   62] noticed, he implicitly defines the measure of an angle as its ratio to a unit angle [cf. De Morgan [1837] pp.     iv, v, 1, 13]. In the same work he provides [pp.    6-13] a defective definition [cf. pp.    6,13] of the radian.

  Hobson[1918, p.   3] defines the measure of an angle as it ratio to a unit angle, without clarifying what is a ratio. In pp.    4,5 of the same work he provides a defective definition of the radian.

  Employing the above-developed theory of ratios, their work may be improved and linked to the measure currently given in the literature [cf. Lebesgue [1966] pp.    38-40] as follows.

  To simplify the subsequent treatment, an equivalence class (of e.g., slss or angles) may be replaced by one of its elements. The intention will be clear from the context. 

   PROPOSITION VI.33 of Euclid [1956, vol. 2, pp.     273-6]. In equal circles
   angles have the same ratio [ratioally equivalent] as the circumferences [arcs] on which they stand,
   whether they stand at the centers or at the circumferences.
   
   There is a (bridgeable) gap in the proof of proposition III.26 of Euclid
   [1956, vol. 2, pp.     56-7] on which the proof of the above proposition is
   based. As a matter of fact, had the proof of proposition III.26 of Euclid
   been completed as its steps (but the last) indicate, it would have proved:
   
   PROPOSITION III.26$^{\prime }$. In equal circles equal angles have equal
   (central) sectors, whether the angles stand at the centers or at the
   circumferences.\newline
   
   Instead of the actually proven proposition, which is:
   
   PROPOSITION III.26 of Euclid [1956, vol. 2, p.    56]. In equal circles equal
   angles stand on equal circumferences [arcs], whether they stand at the
   centers or at the circumferences.
   
   Based on proposition III.26$^{\prime }$, the proof of proposition VI.33 of
   Euclid may be modified to yield:
   
   PROPOSITION VI.33$^{\prime }$. In equal circles angles have the  same ratio [ratioally equivalent]
   as their (central) sectors, whether they stand at the centers or at the
   circumferences.
   
   Accordingly, in measuring angles, ratios  between angles may be replaced by
   ratios  between the corresponding arcs or sectors,  of congruent circles. To continue go through the following steps.

   (i) Let $Ci$ be a circle (the plane region, not just the boundary) of radius $r$, and let $c$ be an arc of $Ci$.
   Consider the set of all (open) polygons inscribed in $c$. To each of these
   polygons there corresponds a sls, the sum of its edges. The set of all of
   these slss is non-empty and is bounded above (see (iv) below), so it has a supremum sls. This
   defines a function $b$ which assigns to each circular arc the corresponding
   supremum sls.
   
   It is easy to see that $b$ is finitely additive, hence for every pair $%
   c_{1},c_{2}$ of arcs of $Ci$, Ratio $(c_1,c_2)=$ Ratio $(b(c_1),b(c_2))$. On the
   other hand finite additivity follows from this equality, by proposition V.24
   of Euclid [1956, vol. 2, p.    183].
   
   Let $a$ be an angle, and let $u$ be a chosen (unit) angle. Regarding $a$ and 
   $u$ as central angles of $Ci$, let $c_{a}$ and $c_{u}$ be, respectively, the
   corresponding $Ci$ arcs, then by proposition VI.33 of Euclid [1956, vol. 2,
   p.    273]: 
   
   \[
   \text{Ratio}\,\,(a,u)= \text{Ratio}\,\, (c_a,c_u)=  \text{Ratio}\,\,(b(c_a),b(c_u))
   \]
   
   Further, let $Ci^{\prime }$ be a circle of radius $r^{\prime }$, and let $%
   c_{a}^{\prime }$ and $c_{u}^{\prime }$ correspond, respectively, to $%
   c_{a}^{{}}$ and $c_{u}^{{}}$. By similarity of triangles and the properties
   of the suprema, $\text{Ratio}\,\, (b(c_u),r)=\text{Ratio}\,\, (b(c'_u),r') $,

  By proposition V.8 of Euclid [1956, vol.2, 149] it is easy to see that if $u$ is the right angle, this ratio is greater than $Ratio(r,r)$ and if $u$ is one half of the right angle, the ratio is less than $Ratio(r,r)$. So, by continuity there is a unique angle (call it ``radian'' and denote it by ``$d$'') which makes this ratio equal to $Ratio(r,r)$, that is :

  $Ratio\left(b\left(c_d\right),r\right)=Ratio(r,r)$, hence, $b\left(c_d\right)=r$. Consequently:
  \begin{eqnarray*}
Ratio\left(a,d\right) &=& Ratio\left(b(c_a\right),b\left(c_d\right)) \notag \\
 &=& Ratio(b(c_a),r) \notag \\
 &=& Ratio(b({c'}_a),r').
\end{eqnarray*} 
 
 (ii) Enrich the vocabulary of (i) by letting $s_{a}$ and $s_{u}$ ($%
 s_{a}^{\prime }$ and $s_{u}^{\prime }$) be, respectively, the sectors
 corresponding to $a$ and $u$ in $Ci$ $(Ci^{\prime })$. Also, let $q$ and $%
 q^{\prime }$ be, respectively, the squares on $r$ and $r^{\prime }$.
 
 By proposition VI.33$^{\backprime }$ [see above]: 
 \[
 Ratio\,\,(a,u)= Ratio \,\,(s_a,s_u) 
 \]
 
 A corollary of this and proposition XII.2 of Euclid [1956, vol. 3, p.    371
 and the appendix of this article] is: 
 \[
  Ratio \,\,(s_u,s'_u)= Ratio \,\,(q,q'),\]
 {hence}
 	\[	Ratio \,\,(s_u,q)= Ratio \,\,(s'_u,q').
 	 \]
 
 In the proof, proposition V.16 of Euclid [1956, vol. 2, p.    164] is made use
 of.

  By common notion 5 [Euclid [1956] vol.1, p.   155] and proposition V.8 of Euclid [1956, vol.2, p.   149] it is easy to see that if $u$ is the right angle, this ratio is less than $Ratio(q,q)$ and if $u$ is double the right angle, the ratio is greater than $Ratio(q,q)$. So, by continuity there is a unique angle (denote it by ``$e$'') which makes this ratio equal to $Ratio(q,q)$, that is:
  
\[Ratio\left(s_e,q\right)=Ratio(q,q),\] 

Hence, $s_e=q$. Consequently:
\begin{eqnarray*}
Ratio\left(a,e\right) &=& Ratio\left(s_a,s_e\right) \notag \\ 
&=& Ratio\left(s_a,q\right) \notag \\ 
&=& Ratio\left(s'_a,q'\right). 
\end{eqnarray*}
(iii) So, the(direct) measure of a in terms of unit angle $d(e)$, to be denoted by ``$m(a)$'' (``$\mu(a)$''), is given by: 
\begin{eqnarray*}
m\left(a\right) &=& Ratio\left(a,d\right)=Ratio(b(c_a),r) \notag \\
\mu\left(a\right) &=& Ratio\left(a,e\right)=Ratio(s_a,q).
\end{eqnarray*}
This defines two injective functions $m$ and $\mu$ from angles to ratios (real numbers). In modern terminology it may be written:

\[a=m\left(a\right)d=\mu(a)e.\] 

Notice that these equations involve no multiplication; $a,d$ and $e$ are angles, not numbers.

   (iv) To figure out the relationship between $m,d$ and, respectively, $\mu ,e$,
   let $p_{i}(p_{o})$ be the sls formed by adding the edges of an inscribed (a
   circumscribed) polygon in (around) $c_{a}$, and let $Rl(x,y)$ be the
   rectangle whose sides are the slss $x$ and $y$:
   \begin{enumerate}
   	\item[(1)] $p_{i}<p_{o}.$
   	\item[(2)] $2s_{a}<Rl(r,p_{o})$. By exhaustion, the r.h.s. may be made arbitrarily
   	close to the l.h.s.
   	\item[(3)] $Rl(r,p_{i})<2s_{a}$. The proof is by contradiction making use of (1)
   	and (2). Again by exhaustion, the l.h.s. may be made arbitrarily close to
   	the r.h.s.
   	\item[(4)] $Rl(r,b(c_{a}))$ $\zeta $ $2s_{a}$. Recall that $\zeta $ is defined in
   	part (v) of II.1 above.
   	\item[(5)]
 Notice that for positive integers $n,n'$ and magnitudes $x,y$ (of the same or different kinds), $Ratio\left(nx,n'x\right)=Ratio(ny,n'y)$. As usual, this common ratio is denoted by ``$n/n'$'' a just ``$n$'' if $n'=1$. Also, for magnitudes $x,x'$ of the same kind, we write 
 $$ ``n/n'Ratio(x,x')'' \,  \text{for} \,  ``n/n'\times Ratio(x,x')''. \text{ \, Accordingly:} $$
 \vspace{-30pt}  
 \begin{eqnarray*} 
n/n'Ratio(x,x') &=& Ratio(nx,n'x)\times Ratio(x,x') \notag \\ 
&=& Ratio(nx,n'x)\times Ratio(n'x,n'x') \notag \\ 
&=& Ratio(nx,n'x')
\end{eqnarray*}
\item[(6)]
\vspace{-20pt}
\begin{eqnarray*}
2\mu \left(a\right) &=& Ratio\left(2s_a,q\right)=Ratio(Rl\left(r,b\left(c_a\right)\right),q)  \notag \\
&=& Ratio\left(b\left(c_a\right),r\right)=m(a)
\end{eqnarray*} 
The equality before the last follows from proposition VI.1 of Euclid [1956, vol.2, P191].
\item[(7)]  With the obvious definition, $m=2\mu $.
\item[(8)]
\begin{eqnarray*}
Ratio\left(a,e\right) &=& \mu \left(a\right)=(1/2)m(a)=(1/2)Ratio(a,d) \notag \\ 
&=& Ratio(a,2d).
 \end{eqnarray*}
So, by proposition V.9 of Euclid [1956, vol.2, p.   153], $e=2d$.
\end{enumerate}
This is the reason behind the factor 2 in [Hardy [1967] p.    317] and [Morrey
[1962] p.    214]. There the computation is correct, while the treatment is
conceptually completely inappropriate [cf. II.1 above, II.5 below, and
Lebesgue [1966] pp.     63-4].

\vskip0.75cm

(v) The above discussion raises the question: Why did not Euclid prove
	 	propositions along the lines of III.26$^{\prime }$ and VI.33$^{\prime }$,
	 	though he has developed all the necessary infrastructure for proving them?
	 	
	 	This question becomes more pressing when it is remembered that   Euclid
	 	[1956] proved:
	 	
	 	PROPOSITION XII.2 of Euclid [1956, vol. 3, pp.     371-8]. Circles are to one
	 	another as the squares on their diameters.

	 	There is a bridgeable gap in Euclid's proof [see the appendix of this
	 	article].

	 An immediate corollary of this proposition is:
	
	 {PROPOSITION X11.}$\boldsymbol{2}\boldsymbol{'}$\textbf{. }Circles have the same ratio to the (respective) squares on their radii.
	
	  Denote this ratio by ``$\pi $'' and denote the right angle by ``$\rho$'', then:
	
	  $\pi =Ratio\left(Ci,q\right)=Ratio(4s_\rho,q)$, hence:
	\[m\left(\rho\right)=2\mu \left(\rho\right)=2Ratio\left(s_\rho,q\right)=(1/2)\pi .\] 
	That is, the measure of the right angle in radians is $(1/2)\pi $. 
	
	\vskip0.75cm
	
	  \textbf{II.5. From Geometry to Analysis.}
	Trigonometric Functions have a long history. Hobson [1918, p.   18] says: ``Until recently, the circular [trigonometric] functions of an angle were defined, not as ratios, but as lengths [slss] having reference to [the corresponding intercepted] arcs of a circle of specified size.'' Also [cf. De Morgan [1837] p.   24]. After providing some details, showing -among other things- that what is meant by ``lengths'' are slss, Hobson [1918, p.   18] continues: ``The advantage of the present mode of definition of the functions as ratios, is that they are independent of the radius of any circle, and are therefore functions of an angular magnitude only.'' For more derails about the ancient Greeks' treatment of the subject, Heath [1963, pp.    273-4, 393-414] may be consulted.
	
	  Though the notion of a ratio was problematic as De Morgan [1836, p.   62] noticed, he [De Morgan [1837] pp.    17-8] and Hobson [1918, pp.    15-6, 18] defined the geometric trigonometric functions as ratios. Improving upon their work by employing the above-developed theory of ratios, it is easy to see that there is a function, to be denoted by ``$Sin'$'', from the set $Ac'$, of all non-zero acute angles, to ratios, defined -making use of the similarity of triangles-  by:
	
	  $Sin'\theta =Ratio$ (opposite side, hypotenuse)
	
	  Where $\theta $ is one of the two acute angles of a right-angled triangle. 
	  
	 % \newpage
	  
	  \begin{figure}[h!]
	  \centering
	  \includegraphics[scale=0.3]{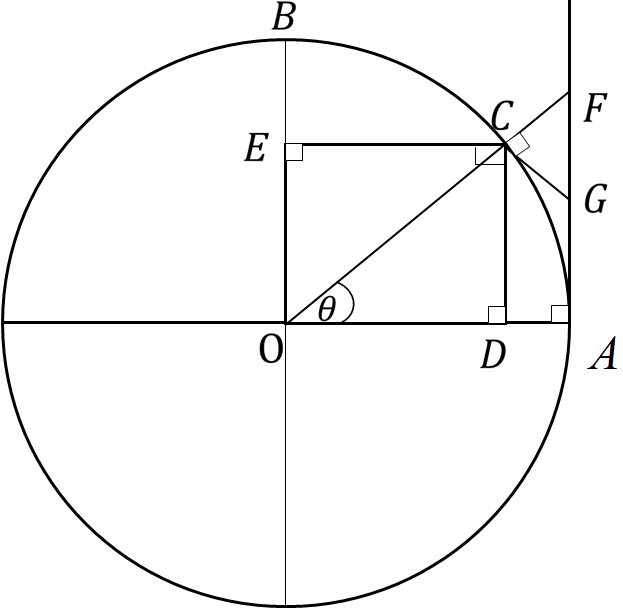}
	  \caption{}
	  \end{figure}

	  Regarding continuity, there is a bijection between $Ac'$ and the points of the arc $\widearc{AB}$ (excluding the points A, B) mapping the angle $\theta $ on the point C as shown in figure1. Also, there is a bijection between these points and the points of the sls $ \overline{OB}$ (excluding the points O, B) mapping the point C on the point E as shown in the same figure. So there is a bijection between $Ac'$ and the set $\left\{Ratio\left(\overline{OE},\overline{OA}\right):E\in \overline{OB}-\left\{O,B\right\}\right\}$.
From this it follows that the $Sin'$ function is a bijection from $Ac'$ onto ]0,1[. Here we feel free to switch back and forth between the language of ratios and the current language of real numbers, as we have shown that the ratios are the positive elements of the E-D field, which is the (geometric) field of the real numbers.
 	
	  Each of the transitions of De Morgan [1837] and Hobson [1918] from the geometric trigonometric functions to their analytic counterparts is again defective. De Morgan [1837, p.   19] considers the analytic trigonometric functions to be abbreviations of their geometric counterparts, while Hobson [1918, pp.    17, 124] confuses the former with the latter. 
	
	  Moreover, each of the treatments of De Morgan [1837, pp.    28-9] and Hobson [1918, pp.    124-5] of the famous ${\mathop{\mathrm{lim}}_{x\to 0} sinx/x=1\ }$, is inadequate. Meanwhile, Hobson's treatment has some merits which we shall build upon in the forthcoming process of employing the above-developed theory of ratios to improve the works of De Morgan and Hobson. 
	
	  (i) For every positive real number $x$, there is an acute angle $\theta \epsilon Ac'$ such that $m\left(\theta \right)<x$. For, regarding figure 1 [cf. Hobson [1918] p.   124]:
	
	  (1) There is a sls $l$ such that $x=Ratio(l,\overline{OA})$.
	
	(2) let F be one of the points on the tangent to the circle at A, such that $\overline{AF}$ is congruent to $l$.
	
	(3) let C be the point of intersection of $\overline{OF}$ and the circumference of the circle.
	
	(4) let G be the point of intersection of the tangent to the circle at C and $\overline{AF}$.
	
	(5) let $\theta $ be the angle $<<F,O,A>,G>$, then:
	\[b\left(c_{\theta }\right)<\overline{CG}+\overline{GA}<\overline{AF}.\] 
	
	(6) 
		$m\left(\theta \right)=Ratio\left(b\left(c_{\theta }\right),\overline{OA}\right)<Ratio\left(\overline{AF},\overline{OA}\right)=Ratio\left(l,\overline{OA}\right)=x. 
	$
	 \vskip0.5cm
	 
	(ii) let $m'$ be the restriction of the measure function $m$ to $Ac'$. Since $m$ is injective, $m'$ is injective too. By (i) and continuity, $m'$ is a bijection from $Ac'$ onto ]$0,\pi /2$[.

	\begin{figure}[htp!]
\centering
\includegraphics[scale=0.4]{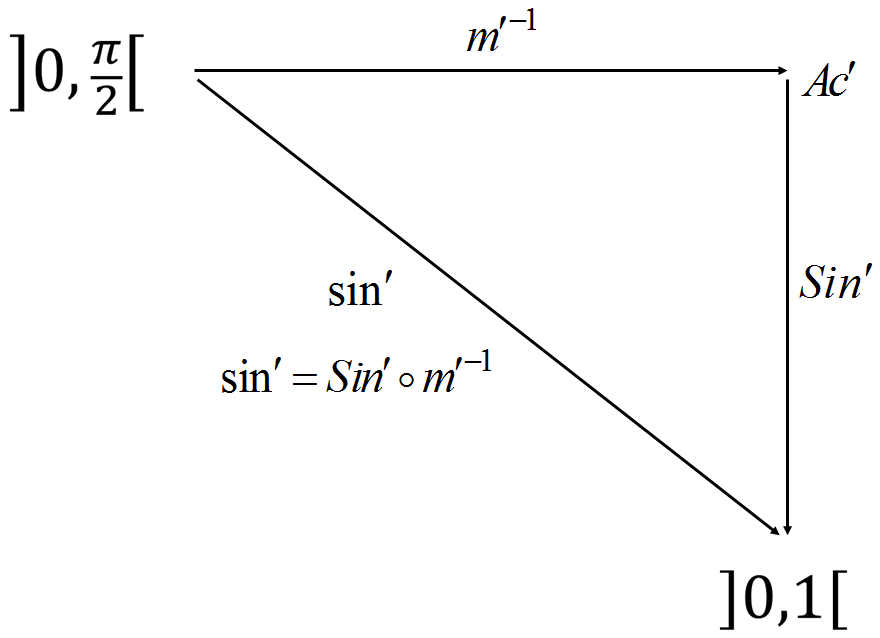}
\caption{}
	\end{figure}
	
	  (iii) Regarding figure 2, it is easy to see that $sin'$ is a bijection from \\  ]\textit{ }$0,\pi /2$ [ onto]\textit{ }$0,1$ [. It is the restriction of the analytic sine function on ]\textit{ }$0,\pi /2$ [.
	
	  The full geometric and analytic trigonometric functions may now be easily defined. 
	\vskip0.5cm
	
	  (iv) To see that $sin'x/x$ tends to 1 as $x$ decreases to 0, let $x\in ]0,\pi /2[$, then there is $\theta \epsilon Ac'$ such that $m\left(\theta \right)=x$, let it be the $\theta $ shown in figure 1. Then,
	  \begin{eqnarray*}
	  sin'x &=& Sin'\theta = Ratio(\overline{CD},\overline{OA}), \notag \\ 
	  x &=& m\left(\theta \right) = Ratio(b\left(c_{\theta }\right),\overline{OA}).
\end{eqnarray*} 
	Hence,
	\begin{eqnarray*}
	sin'x/x &=& Ratio(\overline{CD},\overline{OA})\times Ratio(\overline{OA},b\left(c_{\theta }\right)) \notag \\
	&=& Ratio\left(\overline{CD},b\left(c_{\theta }\right)\right)<Ratio\left(\overline{CA},b\left(c_{\theta }\right)\right)<1.
		\end{eqnarray*} 
	Also, $sin'x/x=Ratio(\overline{CD},b\left(c_{\theta }\right))>Ratio(\overline{CD},\overline{AF})
	=Ratio(\overline{OD},\overline{OA}).$
	\vspace{10pt}
	
	\noindent But $\overline{DA}<\overline{CA}<b\left(c_{\theta }\right)$ and for $\theta <d$ (the radian), $b\left(c_{\theta }\right)<\overline{OA}$. 	
	So, for $\theta <d$,
	\begin{eqnarray*} sin'x/x>Ratio(\overline{OD},\overline{OA}) &=& Ratio(\overline{OA}-\overline{DA},\overline{OA})
	>Ratio\left(\overline{OA}-\overline{CA},\overline{OA}\right) \notag \\
	&>& Ratio(\overline{OA}-b\left(c_{\theta }\right),\overline{OA}) \notag \\ 
	&=& Ratio\left(\overline{OA},\overline{OA}\right)-Ratio(b\left(c_{\theta }\right),\overline{OA}) \notag \\
	&=& 1-m\left(\theta \right)=1-x.
	\end{eqnarray*} 
	Consequently, $1-x<sin'x/x<1$, hence the result.
	
	  Above, we made use of each of ``$<$'' and ``$-$'' in more than one sense, the context clarifies the intention. 
	
	  From this point we may further proceed as usual. 
	\vskip0.5cm
	
	  \textbf{II.6. Alternative Approaches.}
	
	\begin{minipage}{\linewidth}
		\centering
		\begin{minipage}{0.4\linewidth}
			\begin{figure}[H]
				\includegraphics[width=\linewidth]{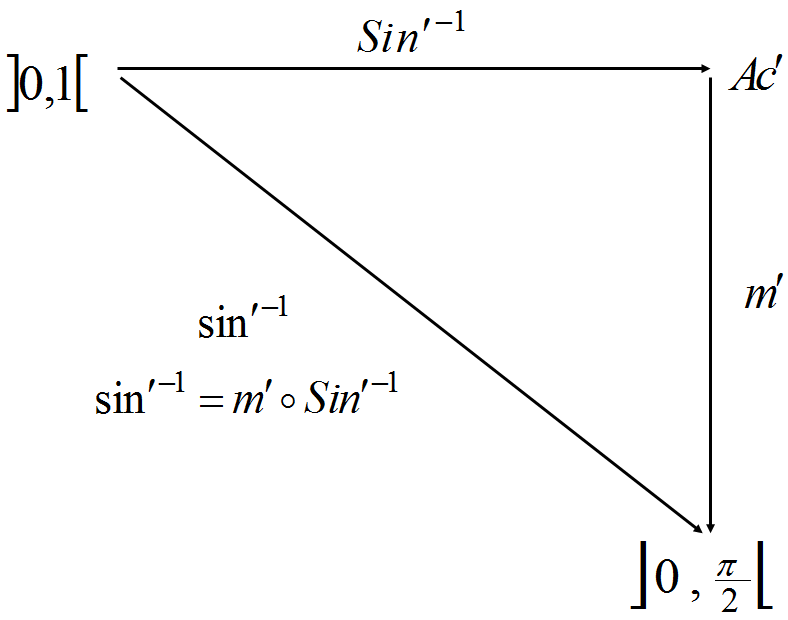}
				\caption{}
			\end{figure}
		\end{minipage}
		\hspace{0.05\linewidth}
		\begin{minipage}{0.4\linewidth}
			\begin{figure}[H]
				\includegraphics[width=\linewidth]{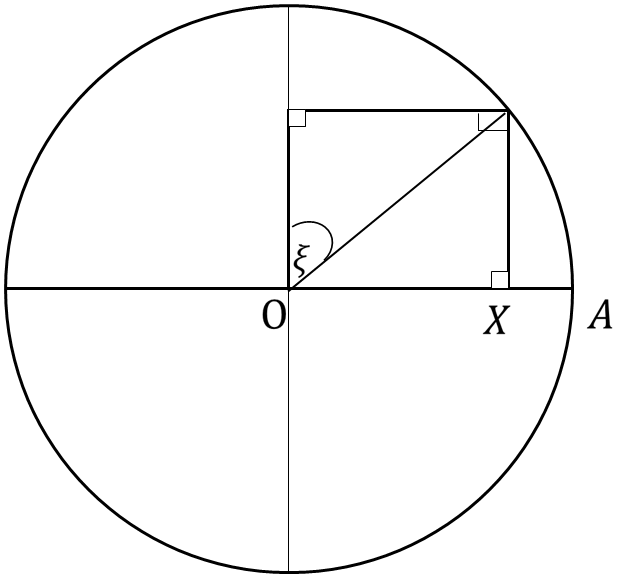}
				\caption{}
			\end{figure}
		\end{minipage}
	\end{minipage}

	  (i) Regarding figures 3,4, let $\overline{OA}$ be the unit $sls$ and let $x\in ]0,1[$. Then there is a point X such that $Ratio\left(\overline{OX},\overline{OA}\right)=x$, hence, ${Sin'}^{-1}\left(x\right)=\xi$, and ${sin'}^{-1}\left(x\right)=m'\left({Sin'}^{-1}\left(x\right)\right)=m'\left(\xi\right)=Ratio(b\left(c_{\xi }\right),\overline{OA})$. 
	
	  So, as is well known, ${sin'}^{-1}\left(x\right)=\int^x_0{dt/\sqrt{1-t^2}}$. From this point we may further proceed as usual [cf. Amer [2017]]. 
	
	  (ii) A great advantage of the Eudoxean theory of direct measure is that it provides a unified method to measure every magnitude by a magnitude of the same kind. In particular, slss are to be measured by slss, angles by angles, circular arcs by circular arcs (of congruent circles), and circular sectors by circular sectors (of congruent circles). For these particular kinds there is another unified method to do the job, namely that expounded by Lebesgue [1966, pp.    19-20, 35-9] and explained (for slss and angles) in I.5 above.
	
	  The Eudoxean direct measures are in fact special cases of ratios, obtained by fixing a magnitude to be the common consequent, or the measuring unit.     
	
	  On the other hand, ratios may be obtained from direct measures by allowing the measuring unit to be varied. For Lebesgue, it is the other way round. Measures are defined first, then ratios are obtained by generalization. Of course, we may return back to measures by specification. So, in both cases it easy to switch back and forth between ratios and direct measures.
	
	  The above definition of the measures function $m$ made use only of the aforementioned kinds. So, this function may be redefined in terms of Lebesguean ratios (with some theory of the real numbers in the background).
	
	  In moving from geometry to analysis we made use of the definition of $\pi $ given above, which needs dealing with kinds other than those mentioned above. But this is dispensable as the equation $m\left(\rho\right)=(1/2)\pi $, or rather, $\pi =2m\left(\rho\right)$, may be taken as the definition of $\pi $. So, the above transition from geometry to analysis may be redone in terms of Lebesguean ratios instead.
	
	  The above definition of the measure function $\mu $ made use also of Eudoxean ratios between circles, circular sectors, and rectangles. Lebesgue [1966, pp.    42-5, 62-3] defined direct measures (areas, in his terminology) of these magnitudes taking a square as the unit magnitude, but making use of a method different from his unified method mentioned above. This different method makes switching back and forth between measures and ratios problematic, if it is possible at all. As it is not clear whether this method in applicable if a circle, for example, is taken as the measuring unit.
	
	  To get around this difficulty, ratios between magnitudes may be replaced by ratios between their respective measures (which are real numbers). This is reasonable, as it gives the same result for Eudoxean ratios. Thus, the term ``ratio'' will have two different senses. The only kind to which they both apply is that of circular sectors of congruent circles. For this kind it may be shown that they both give the same (numerical) result. So, we may freely make use of this term without risking any computational confusion.
	
	  Consequently, all what is done above concerning measuring angles and moving from geometry to analysis may be redone replacing Eudoxean ratios by their Lebesguean counterparts.
	
	  The different approaches to the geometric and the analytic trigonometric functions expounded above are not conceptually equivalent, meanwhile, it may be proved that they are all computationally equivalent.
	
	%\newpage
	\vskip0.75cm
	  \textbf{III. CONCLUDING REMARKS}

	\textit{Retracing Elementary Mathematics} [Henkin et al. 1962] is the
	outcome of three Summer Institutes for Teachers of Mathematics, organized by
	the National Science Foundation, and taught and directed by the authors
	[Henkin et al. [1962] p.    vii]. It is supposed to initiate a new beginning,
	and is proposed to be ``Used as a college text [...] for a one-year course
	of the type commonly given under a title such as ``Foundations of
	Mathematics,'' or ``\textit{Fundamental Concepts of Mathematics}.''
	[emphasis added].''  [Henkin et al. [1962] p.    viii]. 
	
	Nevertheless the book deals essentially with Arithmetization of the real numbers from ``the modern deductive point of view''
	[Henkin et al. [1962] p.    viii] considering the \textit{Fundamental Concepts of
		Elementary Mathematics} such as the \textit{length of a sls} and the\textit{ measure of an
		angle } to be already known, with no word concerning retracing them [cf. Henkin et
	al. [1962] pp.     316, 328, 332].
	
	This is compatible with the current practice in mathematics. It is well
	known that treatises on the structure of the real numbers go along the lines of
	Henkin et al [1962] with varying degrees of formality; treatments of
	Euclidean geometry are based on the real numbers or on Hilbert [1950, pp.     23-36]'s theory of ratio and
	proportion, or variations thereof [cf. Borceux [2014, pp.     vii-ix, 305,
	351,352], Hartshorne [2000, pp.     2,3,168], Meyer [2006, pp.     23,29], Sossinsky
	[2012, pp.     1,2], Szmielew [1983, pp.     72-84, 109, 180, 181], Tarski [1959], Tarski and
	Givant [1999, p.    211] and Venema [2006, pp.     ix, 398-410]]; dealing with
	measure (even of slss) is well known to presuppose the real numbers; ... . The
	result is that -to the best of my knowledge- current mathematics almost
	entirely ignores Eudoxus' theory of ratio and proportion; the few exceptions which deal
	with it, treat it as if it were merely a piece of antique. This being so,
	though,
			beside playing an essential role in mathematics, this theory is the medium through which the relationships among different kinds may be studied, it is the medium which makes theoretical natural since possible. For example, it is wrong to say that density is mass divided by volume. Mass and volume are magnitudes (of different kinds), how would one of them be divided by the other? In fact, density is the (Eudexean) ratio of a mass to a unit mass divided by the (Eudexean) ratio of a volume to a unit volume. Ratios (unlike magnitudes) may be divided by one another. 
	
	  Hilbert [1950, pp.    23-36]'s theory of ratio and proportion is obviously insufficient. It did not help Hilbert [1950] to prove the usual formula for the area of a triangle, so, this formula was taken [Hilbert [1950] p.   41] as a definition, without referring to Eudoxus' theory in which it may be proved. 
	
	  On the other hand there may be some reservations against Eudoxus' theory of ratio and proportion as it makes use of the natural numbers, which are considered to be \textit{alien }to geometry. To get around this difficulty, we may follow Hilbert [1950, pp.    29-32] to define a \textit{pure }geometric ordered field. If, in addition, we adopt the second order axiom of continuity (which Hilbert [1950, P36] grants that it is valid in the ordinary analytic geometry of space), this field will be isomorphic to the field of real numbers. With second order logic, the natural numbers may be defined as a subset of this field, so, in a sense, they will be geometrized. Making use of them, the Eudexean ratios, hence the E-D field, may be (geometrically) defined.
	
	 Thus, every thing will be nicely settled, paving the way for incorporating the Eudoxean theory of ratio and proportion in contemporary mathematics. Along with this, mathematics,  its history, philosophy and pedagogy should be radically reconsidered.   
	\newpage
	
		\textbf{ACKNOWLEDGMENTS }
		
		This article has been
		developed through several seminar talks in the Department of Mathematics,
		Faculty of Science, Cairo University. It is my pleasure to express my deep
		gratitude to the friends who invited me to give talks in the seminars they
		are leading, Professors: Annaby, Ghaleb, Megahed, and Youssef.
		
	 \vskip0.5cm

	{\large APPENDIX}

	The gap in each of the proofs of proposition XII.2 of Euclid [1956, vol. 3,
	pp.     371-8; cf. II.4 above] attributed to Eudoxus and Legendre, is the
	unjustified assumption of the existence of the fourth proportional. This
	error is recently reiterated  in the proofs of a version of this proposition
	[Borceux [2014] pp.     35-7, 100-2]. Concerning this existence, the arguments
	attributed to De Morgan [see Euclid [1956] vol. 2, p.    171] and Simson [see
	Euclid [1956] vol. 3, p.    375] do not solve the problem. Nor would
	proposition VI.12 of Euclid [1956] vol. 2, p.    215] help, for it proves the
	existence of the fourth proportional only for straight line-segments.
	
	Following is a proof of proposition XII.2 of Euclid which does not make use
	of this problematic assumption. It is so simple that it may replace the
	proof given by Euclid in his book.
	
	Let $c_{1}$ and $c_{2}$ be circles, and let $s_{1}$ and $s_{2}$ be the
	squares on their respective diameters.
	
	To prove that $<c_{1},c_{2}>=_{E}<s_{1},s_{2}>$, assume the contrary. Then
	by definition V.5 of Euclid [1956, vol. 2, p.    114], there are positive
	integers $n_{1}$ and $n_{2}$ such that:
	\begin{description}
		\item  (i) $n_{1}c_{2}>n_{2}c_{1}$ \quad and $\quad n_{1}s_{2}\leq
		n_{2}s_{1} $,
		\item  (ii) $n_{2}c_{1}>n_{1}c_{2}$ \quad and $\quad n_{2}s_{1}\leq
		n_{1}s_{2}$,
		\item  (iii) $n_{1}c_{2}=n_{2}c_{1}$ \quad and $\quad n_{1}s_{2}<n_{2}s_{1}$,
		\item  or
		\item  (iv) $n_{2}c_{1}=n_{1}c_{2}$ \quad and $\quad n_{2}s_{1}<n_{1}s_{2}$.
	\end{description}
	
	Since (ii) and (iv) may, respectively, be obtained from (i) and (iii) by
	interchanging the indices 1 and 2, it suffices to deal with (i) and (iii).
	
	Let $p_{1}$ and $p_{2}$ be similar polygons inscribed in $c_{1}$ and $c_{2}$
	respectively. By proposition XII.1 of Euclid [1956, vol. 3, p.    369]: 
	\begin{equation}
	<p_{1},p_{2}>=_{E}<s_{1},s_{2}>  \tag{*}
	\end{equation}
	
	Assume (i), then by (*) and definition V.5 of Euclid, 
	\[
	n_{1}p_{2}\leq n_{2}p_{1}<n_{2}c_{1}<n_{1}c_{2} 
	\]
	hence, 
	\[
	(n_{1}c_{2}\stackrel{.}{-}n_{2}c_{1})<(n_{1}c_{2}\stackrel{.}{-}n_{1}p_{2}), 
	\]
	which entails a contradiction. For, in view of exhaustion, the right-hand side
	may be made less than the left-hand side, by choosing an appropriate $p_{2}$.%
	\newline
	
	To deal with the other case, let $p_{1}^{\prime }$ and $p_{2}^{\prime }$ be
	similar polygons inscribed in $c_{1}$ and $c_{2}$ respectively. As above,
	\[
	<p_{1}^{\prime },p_{2}^{\prime}>=_{E}<s_{1},s_{2}> 
	\]
	so by (*) 
	\[
	<p_{1}^{\prime },p_{2}^{\prime}>=_{E}<p_{1},p_{2}> 
	\]
	
	Moreover, let $p_{1}^{\prime },p_{2}^{\prime }$ properly include $%
	p_{1}^{{}},p_{2}^{{}}$ respectively, then by propositions V.16, 17 of Euclid
	[1956, vol. 2, pp.     164-6], 
	\[
	<(p_{1}^{\prime }\stackrel{.}{-}p_{1}),(p_{2}^{\prime}\stackrel{.}{-}p_{2})>%
	=_{E}<p_{1},p_{2}> 
	\]
	
	Finally, assume (iii), then, 
	\[
	n_{2}p_{1}\stackrel{.}{-}n_{1}p_{2}<n_{2}p_{1}^{\prime }\stackrel{.}{-}%
	n_{1}p_{2}^{\prime }<n_{1}c_{2}\stackrel{.}{-}n_{1}p_{2}^{\prime}, 
	\]
	which entails a contradiction.
	
	This completes the proof.
	
	\vskip1.0cm
	
	 \textbf{\Large BIBLIOGRAPHY}

    \begin{description} 
    	\item{}  Al-Khwarizmi, Muhammad Ibn Musca [1939]; (Hisab) \textit{Algabr
    		wa Al-Mukabala} [(Calculus) of Restoration and Cancellation]; edited with
    	introduction and commentary by Mosharrafa, A.M. and Ahmed, M.M.;
    	Publications of the Faculty of Science, the Egyptian (now Cairo) University,
    	no. 2; (In Arabic).
    	\item{}  Amer, Mohamed A. [1980]; \textit{The Relationship of Mathematics
    		to Reality as a Basis of the Philosophy of Mathematics}; Scientia, An.
    	LXXIV, vol. 115, issue 9-12, pp.     563-578.
    	\item{} Amer, Mohamed A. [2017]; \textit{What is the "$x$" Which Occurs in "$sin x$"?} Available Online at  arxiv>math>arxiv:1610.07867V2.
    	
    	\item{} Amer, Mohamed A. [2022]; \textit{Aristotelian Assertoric Syllogistic}; Springer Briefs in Philosophy; Springer Nature Switzerland AG.
    	\item{}  Anton, Howard; Bivens, Irl C.; and Davis, Stephen [2009]; 
    	\textit{Calculus, Single Variable}; 9th ed.; John Wiley \& Sons, Inc..
    	\item{}  Avigad, Jeremy; Dean, Edward; and Mumma, John [2009]; \textit{A
    		Formal System for Euclid's Elements}; The Review of Symbolic Logic; vol. 2,
    	no. 4, pp.     700-768.
    	\item{}  Baldwin, John T. [2013]; \textit{Formalization, Primitive
    		Concepts, and Purity}; The Review of Symbolic Logic, vol. 6, no. 1, pp.    
    	87-128.
    	\item{}  Bartle, Robert G. and Tulcea, C. Ionescu [1968]; \textit{Calculus%
    	}; Scott, Foresman and Company.
    	\item{}  Boole, George [1948]; \textit{The Mathematical Analysis of Logic}%
    	; Oxford, Basil Blackwell.
    	\item{}  Borceux, Francis [2014]; \textit{An Axiomatic Approach to
    		Geometry, Geometric Trilogy I}; Springer.
    	\item{}  Borsuk, Karol and Szmielew, Wanda [1960]; \textit{Foundations of
    		Geometry}; North-Holland Publishing Company.
    	\item{}  Bos, Henk J. M. [1997]; \textit{Lectures in the History of
    		Mathematics}; AMS History of Mathematics; vol. 7.
    	\item{}  Buck, R. Creighton [1965]; \textit{Advanced Calculus}; 2nd ed.;
    	McGraw-Hill Book Company.
    	\item{}  Davis, Martin [2003]; \textit{Exponential and Trigonometric
    		Functions- From the Book}; The Mathematical Intelligencer; vol. 25, no. 1,
    	pp.     5-7.
    	
    	\item{} Dedekind, Riched [1901]; \textit{Continuity and Irrational Numbers}, part I of \textit{Essays on the Theory of Numbers}; translated by Beman, Wooster Woodruff; The Open Court Publishing Company. Also available Online as a Project Gutenberg [EBook\#21016], 2007.
    	\item{} De Morgan, Augustus [1836]; \textit{The connexion of Number and Magnitude}; London: Printed for Taylor and Walton, Booksellers and Publishers to the University of London; Available Online at  http://www.archive.org/details/cnnexionofnumbe00demorich; Digitized by Internet Archive in 2007 with funding from Microsoft Corporation.
    	\item{} De Morgan, Augustus [1837]; \textit{Elements of Trigonometry and Trigonometrical Analysis}; London: Printed for Taylor and Walton, Booksellers and Publishers to the University of London; Available Online via Google Book Search at http://books.google.com/.
    	    	\item{}  Descartes, Ren\textit{\'{e}} [1954]; \textit{The Geometry of
    		Ren\'{e} Descartes}; translated from the French and Latin by David Eugene
    	Smith and Marcia L. Latham; Dover Publications, Inc..
    	\item{}  Eberlein, W. F. [1966]; \textit{The Circular Function(s)};
    	Mathematics Magazine; vol. 39, no. 4, pp.     197-201.
    	\item{}  Euclid [1956]; \textit{The Thirteen Books of Euclid's Elements};
    	translated with introduction and commentary by sir Thomas L. Heath; 2nd ed.
    	unabridged; Dover Publications, Inc. (3 Volumes).
    	\item{}  Euler [2000]; \textit{Foundations of Differential Calculus};
    	translated by John D. Blanton; Springer.
    	\item{}  Feferman, Solomon [1999]; \textit{Does Mathematics Need New
    		Axioms?}; American Mathematical Monthly; vol. 106, no. 2, pp.     99-111.
    	\item{}  Feferman, Solomon; Friedman, Harvey M.; Maddy, Penelope; and
    	Steel, John R. [2000]; Does \textit{Mathematics Need New Axioms?}; The
    	Bulletin of Symbolic Logic; vol. 6, no. 4, pp.     401-446.
    	\item{} Frege, Gottlob [1968]; \textit{The Foundations of Arithmetic}; German text with English translation by J. L. Austin; 2nd. ed.; Northwestern University Press.
    	\item{}Frege, Gottlob [2016]; \textit{Basic Laws of Arithmetic}; translated and edited by Philip A. Ebert and Marcus Rossberg, with Crispin Wright; Oxford University Press; Vols. I \& II.
    	\item{}  Gearhart, William B. and Shultz, Harris S. [1990]; \textit{The
    		Function }$\sin x/x$; The College Mathematics Journal; vol. 21, no. 2, pp.    
    	90-99.
    	\item{}  Gillman, Leonard [1991]; $\pi $\textit{\ and the limit of }$%
    	(\sin \alpha )/\alpha $; American Mathematical Monthly; vol. 98, no. 4, pp.    
    	346-349.
    	\item{}  G\"{o}del, Kurt [1964a]; \textit{Russell's Mathematical Logic};
    	in Benacerraf, Paul and Putnam, Hilary eds.; Philosophy of Mathematics;
    	Prentice-Hall, Inc.; pp.     211-232.
    	\item{}  G\"{o}del, Kurt [1964b]; \textit{What is Cantor's Continuum
    		Problem?}; in Benacerraf, Paul and Putnam, Hilary eds.; Philosophy of
    	Mathematics; Prentice-Hall, Inc.; pp.     258-273.
    	\item{}  Greenberg, Marvin Jay [1980]; \textit{Euclidean and
    		Non-Euclidean Geometries}; 2nd ed.; W.H. Freeman and Company.
    	\item{}  Griffiths, Phillip A. [2000]; \textit{Mathematics at the Turn of
    		the Millennium}; American Mathematical Monthly; vol. 107, no. 1, pp.     1-14.
    	\item{}  Hardy, G. H. [1967]; \textit{A Course of Pure Mathematics}; 10th
    	ed.; Cambridge University Press. (1st ed. 1908).
    	\item{}  Hardy, G. H. [1969]; \textit{A Mathematician's Apology} (Foreword
    	by C.p.    Snow); Cambridge University Press.
    	\item{}  Hartshorne, Robin [2000]; \textit{Geometry: Euclid and Beyond};
    	Springer - Verlag New York, Inc.
    	\item{} Hass, Joel; Heil, Christopher; and Weir, Maurice D. [2020]; \textit{Thomas' Calculus} Based on the original work by George B. Thomas Jr. (-14$^{th}$ed. in SI Units,) Pearson Education, Inc. (Pearson Ed. Limited, UK-  www.pearsonglobaleditions.com)    	
    	\item{}  Heath, Sir Thomas L. [1963]; \textit{A Manual of Greek
    		Mathematics}; Dover Publications, Inc.
    	\item{}  Henkin, Leon; Smith, Norman W.; Varineau, Verne J.; and Walsh,
    	Michael J. [1962]; \textit{Retracing Elementary Mathematics}; The Macmillan
    	Company.
    	\item{}  Hilbert, David [1950]; \textit{Foundations of Geometry};
    	translated by E.J. Townsend; The Open Court Publishing Companry. Also
    	available Online as a Project Gutenberg Ebook, file 17384-pdf.pdf.
    	\item{} Hobson, E. W. [1918]; \textit{A Treatise on Plane Trigonometry}; 4$^{th}$ ed.; Combridege: at the University Press; 1$^{st}$ ed. 1891, 3$^{rd}$ ed. (revised and enlarged) 1911; Also available Online at\\ http://www.archive.org/detalis/treatiseonplanet00hobs;Digitized by the Internet archive in 2008 with funding from Microsoft Corporation. 
    	
    	\item{}  Hobson, E. W. [1957]; \textit{The Theory of Functions of a Real
    		Variable \& the Theory of Fourier's Series}; 3rd ed.; Dover Publications,
    	Inc. (Volume 1).
    	\item{}  Hughes-Hallett, Deborah; Gleason, Andrew M.; Mc Callum, William G.; Flath, Daniel E.; Lomen, David O.; Lovelock, David; and Tecosky-Feldman, Jeff [2005]; \textit{Calculus, Single
    		Variable}; 4th ed.; John Wiley \& sons, Inc.
    	\item{}  Kac, Mark and Ulam, Stanislaw [1971]; \textit{Mathematics and
    		Logic}; a Pelican Book.
    	\item{}  Kneale, William and Kneale, Martha [1966]; \textit{The
    		Development of Logic}; Oxford University Press.
    	\item{}  Lavers, Gregory [2009]; \textit{Benacerraf's Dilemma and
    		Informal Mathematics}; The Review of Symbolic Logic, vol. 2, no. 4, pp.    
    	769-785.
    	\item{}  Lebesgue, Henri [1966]; \textit{Measure and the Integral};
    	translated and edited with a bibliographical essay by Kenneth O. May;
    	Holden-Day, Inc.
    	\item{}  \L ukasiewicz, Jan [1998]; \textit{Aristotle's Syllogistic From
    		the Standpoint of Modern Formal Logic}; 2nd ed. enlarged; Oxford University
    	Press Inc., New York.
    	\item{}  Maddy, Penelope [2008]; \textit{How Applied Mathematics Became
    		Pure}; The Review of Symbolic Logic, vol. 1, no. 1, pp.     16-41.
    	\item{}  Marquis, Jean-Pierre [2013]; \textit{Categorical Foundations
    		of Mathematics}; The Review of Symbolic Logic, vol. 6, no. 1, pp.     51-75.
    	\item{}  Meyer, Walter [2006]; \textit{Geometry ant its Applications};
    	Elsevier Academic Press, 2nd ed.
    	\item{}  Moise, Edwin [1964]; \textit{Elementary Geometry From an
    		Advanced Standpoint}; Addison-Wesley Publishing Company, Inc.
    	\item{}  Morrey, Charles B. Jr. [1962]; \textit{University Calculus with
    		Analytic Geometry}; Addison-Wesley Publishing Company, Inc.
    	\item{}  Moschovakis, Yiannis N. [1995]; \textit{Review of:
    		Computability, by D.S. Bridges, Spr.-Ver. 1994}; American Mathematical
    	Monthly; vol. 102, no. 8, pp.     752-755.
    	\item{}  Protter, Murray H. and Morrey, Charles B. Jr. [1963]; \textit{%
    		Calculus with Analytic Geometry, A First Course}; Addison-Wesley Publishing
    	Company, Inc.
    	\item{}  Richman, Fred [1993]; \textit{A Circular Argument}; The College
    	Mathematics Journal; vol. 24, no. 2, pp.     160-162.
    	\item{}  Rose, David A. [1991]; \textit{The Differentiability of }$\sin x$%
    	; The College Mathematics Journal; Vol. 22, no. 2, pp.     139-142.
    	\item{}  Sher, Gila [2013]; \textit{The Foundational Problem of Logic};
    	The Bulletin of Symbolic Logic; vol. 19, no. 2, pp.     145-198.
    	
    		\item{} Simons, Peter M. [1992]; \textit{Frege's Theory of Real Numbers}; ch. 5 of Simons, Peter; Philosophy and Logic in Central Europe from Bolzano to Tarski; Springer Science+Business Media Dordrecht; pp.     117-141.
    	\item{} Simons, Peter M. [2011]; \textit{Euclid's Context Principle}; Hermathena; no. 191, Winter 2011 [published 2014], pp.     5-24.
    
    	\item{} Sinkevich, Galina I. [2014]; \textit{Concepts of a number of C. Meray, E. Heine, G. Cantor, R. Dedekind and K. Weierstrass}; Technical Transactions, Fundamental Sciences, Krakow, 1-NP, pp. 211-223.
    	\item{} Smith, D. E. [1958]; \textit{History of Mathematics}; Dover Publications, Inc. (Volume II).
    	\item{}  Sossinsky, A. B. [2012]; \textit{Geometries}; AMS student
    	Mathematical Library, Vol. 64.
    	\item{}  Stein, H. [1990]; \textit{Eudoxus and Dedekind: On the Ancient
    		Greek Theory of Ratios and its Relation to Modern Mathematics}; Synthese;
    	vol. 84, no. 2, pp.     163-211.
    	\item{}  Szmielew, Wanda [1983]; \textit{From Affine to Euclidean
    		Geometry, an Axiomatic Approach}; edited, prepared for publication and
    	translated from the Polish by Maria Moszy\'{n}ska; PWN-Polish Scientific
    	Publishers - Warszawa.
    	\item{} Tarski, Alfred [1956]; \textit{The Concept of Truth in Formalized Languages}; ch. VIII of Tarski, Alfred; Logic, Sematics, Metamathematics; Oxford. At the Clarendon Press; pp.     152-278. 
    	\item{}  Tarski, Alfred [1959]; \textit{What is Elementary Geometry?}; in
    	Henkin, L., Suppes, p.    and Tarski, A. eds.; The Axiomatic Method, with
    	Special Reference to Geometry and Physics; North Holland, pp.     16-29.
    	\item{}  Tarski, Alfred and Givant, Steven [1999]; \textit{Tarski's
    		System of Geometry}; The Bulletin of Symbolic logic, vol. 5, no. 2, pp.    
    	175-214.
    	\item{}  Ungar, Peter [1986]; \textit{Reviews} (of three Calculus Books);
    	American Mathematical Monthly; vol. 93, no. 3, pp.     221-230.
    	\item{}  Vafea, Flora K. [1998]; \textit{Topics in Ancient Egyptian
    		Mathematics}; M.Sc. Thesis, Dept. of Math., Faculty of Science, Cairo Univ.,
    	Egypt.\newline
    	\item{}  Venema, Gerard A. [2006]; \textit{The Foundations of Geometry};
    	Pearson Prentice Hall.
    \end{description}

  Department of Mathematics
  
  Faculty of Science
  
  Cairo University
  
  Giza - Egypt
  
  m0amer@hotmail.com
  
  amer@sci.cu.edu.eg

\end{document}